\documentclass[12pt] {article}
\usepackage{amssymb,amsxtra,amsmath,amsthm}
\def \NN{{\mathbb{N}}}
\def \ZZ{{\mathbb{Z}}}
\def \QQ{{\mathbb{Q}}}

\def \FF{{\mathbb{F}}}
\def \PP{{\mathbb{P}}}

\def \CCC{{\mathcal{C}}}
\def \OOO{{\mathcal{O}}}

\def \aa{{\mathfrak a}}
\def \bb{{\mathfrak b}}
\def \cc{{\mathfrak c}}

\def \mm{{\mathfrak m}}
\def \pp{{\mathfrak p}}
\def \qq{{\mathfrak q}}
\def \rr{{\mathfrak r}}
\def \ss{{\mathfrak s}}

\def\mod{\mathop{\mathrm{mod}}}

\def\star{\mathop{\mbox{\Large$*$}}}
\def\Star{\mathop{\mbox{\LARGE$*$}}}

\def\dim{\mathop{\mathrm{dim}}}

\def\Gamma{\varGamma}
\def\Delta{\varDelta}

\begin{document}

\begin{center}
{\Large {\bf The cusp amplitudes and quasi-level of a congruence subgroup of $SL_2$ over any Dedekind domain}}\\
\bigskip
{\tiny {\rm BY}}\\
\bigskip
{\sc A. W. Mason}\\
\bigskip
{\small {\it Department of Mathematics, University of Glasgow\\
Glasgow G12 8QW, Scotland, U.K.\\
e-mail: awm@maths.gla.ac.uk}}\\
\bigskip
{\tiny {\rm AND}}\\
\bigskip
{\sc Andreas Schweizer}\\
\bigskip
{\small {\it Institute of Mathematics, Academia Sinica\\
6F, Astronomy-Mathematics Building\\
No. 1, Sec. 4, Roosevelt Road\\
Taipei 10617, TAIWAN\\
e-mail: schweizer@math.sinica.edu.tw}}
\end{center}
\begin{abstract}
We extend some algebraic properties of the classical modular group 
$SL_2(\ZZ)$ to equivalent groups in the theory of Drinfeld modules,
in particular properties which are important in the theory of modular 
curves.
\par 
We study cusp amplitudes and the level of a (congruence) subgroup
of $SL_2(D)$ for any Dedekind domain $D$, as ideals of $D$. 
In particular, we extend a remarkable result of Larcher.
\par
We introduce finer notions of {\it quasi-amplitude} and 
{\it quasi-level} which are not required to be ideals and
encode more information about the subgroup.
\par
Our results also provide several new necessary conditions for
a subgroup of $SL_2(D)$ to be a congruence subgroup. 
\end{abstract}

\subsection*{Introduction}

\par This paper is part of an ongoing project which aims to extend 
algebraic results from
the classical modular group, $SL_2(\ZZ)$, to equivalent groups
occurring in the context of Drinfeld modular curves [Ge]. We are
especially concerned with those results which have proved applicable
to the classical theory of modular forms. Our hope is that our
results will prove useful to experts in the theory of Drinfeld
modules. Here we are primarily concerned with the {\it cusp
amplitudes} and {\it level} of a subgroup of such a group. One of
our principal aims is to demonstrate that these concepts, as well as
much of the classical theory of congruence subgroups, can be
extended to linear groups defined over {\it any} Dedekind domain
$D$. \par Congruence subgroups, beginning in the $19$th century with
$SL_2(\ZZ)$, are usually defined for a matrix group with entries in
an arithmetic Dedekind domain $A$ (or, more generally, an {\it
order} in such a domain). The definition involves a non-zero
$A$-ideal. Since every proper quotient of $A$ is finite, every
congruence subgroup is necessarily of finite index. This leads
naturally to the converse question (for the particular matrix
group), namely the celebrated {\it Congruence Subgroup Problem},
which has attracted considerable attention for many years. For $D$
which are not arithmetic, i.e. $D$ which might have proper infinite
quotients, there is no longer a close connection between the
congruence subgroups of $SL_2(D)$ and its finite index subgroups.
For example, if $k$ is an infinite field, then $SL_2(k[t])$ has
infinitely many congruence subgroups and no proper finite index
subgroups.
\par
 The groups $SL_2(D)$ include a number of very important special cases,
for example, the classical modular group $SL_2(\ZZ)$, the Bianchi
groups, where $D$ is the ring of integers in an imaginary quadratic
number field, and groups occurring in the context of Drinfeld
modular curves [Ge]. Other examples are $SL_2$ over a ring of
$S$-integers of a number field or of a function field in one
variable over any constant field, $SL_2$ over any local ring and
$SL_2$ over any principal ideal domain.
\par
When generalizing concepts from the classical modular group
$SL_2(\ZZ)$ to $SL_2(R)$ for more general rings $R$, one sometimes
has to make a choice between loss of structure and loss of
information. For example, let $H$ be a subgroup of $SL_2(R)$.
Then the set of all $a\in R$ such the translation matrix
${1\ a\choose 0\ 1}$ is contained in $H$ is an additive subgroup
of $R$. Hence in the special case $R=\ZZ$ this is automatically
an ideal. In the general case one can choose between considering
this set and losing the ideal structure or considering the biggest
ideal contained in this set and thus losing information about $H$.
\par
This is one reason why it is sometimes difficult to extend results
that hold for $SL_2(\ZZ)$ to more general groups $SL_2(R)$. Either the
associated objects do not have enough structure and hence one cannot
prove enough about them, or they have the necessary structure but they
don't encode all the desired information.
This also indicates that one should look at both of the possible
generalizations and how they are related.
\par
Let us make this precise by introducing the most important instance
of this dichotomy structure versus information.
\par
Let $R$ be a commutative ring with identity and let $H$ be a subgroup
of $SL_2(R)$. The classical definition of the {\it level} $l(H)$ of
$H$ opts for structure over information and defines it as the largest
$R$-ideal $\qq$ such that for every $a\in\qq$ the translation matrix
${1\ a\choose 0\ 1}$ is contained in every conjugate of $H$ in
$SL_2(R)$. We emphasize information over structure and call the set
of all $a\in R$ such that the translation matrix ${1\ a\choose 0\ 1}$ 
is contained in every conjugate of $H$ in $SL_2(R)$ the 
{\it quasi-level} of $H$, denoted by $ql(H)$. Implicitly the quasi-level
already occurs for example in the proof of Proposition 1 in [Se1].
\par
Obviously, the level of $H$ is the largest $R$-ideal contained in the
quasi-level of $H$. Clearly for the special case $R=\ZZ$ we always have
$ql(H)=l(H)$. In general however, as we shall see, they can differ by
``as much as possible". By definition the level of a subgroup is not
merely an additive subgroup of $R$ and it can be shown to have much
stronger properties than the quasi-level.
However as a parameter for studying $H$ it is less useful
than $ql(H)$, mainly for the following reason. If $H$ is of finite
index in $SL_2(R)$ then $ql(H)$ is non-zero. On the other hand for
some $R$ (including $R=k[t]$, where $k$ is a finite field) it is
known that $SL_2(R)$ has infinitely many finite index subgroups of
level zero. Moreover for any finite index subgroup $H$ important
information is provided by equations involving the size of the
quasi-level of $H$ and its index in $SL_2(R)$.
\par
Now let $D$ be any {\it Dedekind domain} with quotient field $F$.
The group $SL_2(D)$ acts on the projective line
$\PP^1(F)=F\cup\{\infty\}$ as a set of linear fractional
transformations. Let $H$ be a subgroup of $SL_2(D)$. Clearly $H$
acts on $O_{\infty}=\{g(\infty):g \in SL_2(D)\}$, the
$SL_2(D)$-orbit containing $\infty$. We refer to the orbits of the
$H$-action on $O_{\infty}$ as the $H$-{\it cusps}.
\par
The classical definition of cusp amplitude again prefers structure
over information and defines the {\it cusp amplitude} of the $H$-cusp
containing $g(\infty)$ as a certain $D$-ideal $\cc(H,g)$ associated
to the stabilizer of $g(\infty)$ in $H$. See Section 2 for details.
\par
As the cusp amplitude often loses too much information about the
stabilizer of the cusp, we also introduce what we call the
{\it quasi-amplitude}. This represents the other possibility to
generalize the notion of cusp-amplitude from $SL_2(\ZZ)$ to
$SL_2(D)$, going for more information and less structure by dropping
the requirement that it is an ideal and considering all translation
matrices in $H^g=\{g^{-1}hg\ :\ h\in H\}$.
\par
Again, in general the quasi-amplitude $\bb(H,g)$ is only an additive
subgroup of $D$. Obviously, the cusp amplitude $\cc(H,g)$ is the
largest ideal contained in $\bb(H,g)$. In situations where one can
control how far the two can be apart, one then obtains results that
generalize theorems for the classical modular group.
\par
Let $\mathcal{A}(H)$ denote the set of all the cusp amplitudes of
$H$. It is easy to see that the intersection of all these cusp
amplitudes is nothing else than the level of $H$. In our first
principal result we prove that for one important class of subgroups
the set $\mathcal{A}(H)$ has the following surprising properties.
\\ \\
{\bf Theorem A.} \it
Let $H$ be a subgroup of $SL_2(D)$ and let
$$\cc_{\mathrm{min}}=\displaystyle{\bigcap_{\qq \in \mathcal{A}(H)}}\qq
\ \ \ \ \mbox{\it and}\ \ \ \
\cc_{\mathrm{max}}=\displaystyle{\sum_{\qq \in \mathcal{A}(H)}}\qq.$$
\noindent If $H$ is a congruence subgroup, then
$\cc_{\mathrm{min}},\;\cc_{\mathrm{max}}\in \mathcal{A}(H).$
\par
In particular, if $H$ is a congruence subgroup, there is an $H$-cusp
in the $SL_2(D)$-orbit of $\infty$ whose cusp amplitude equals the
level of $H$.
\rm
\\ \\
This extends the results of a remarkable paper of Larcher [La] for the
particular case $D=\ZZ$. For this case Stothers [St] has also proved
that $\mathcal{A}(H)$ has a minimum by an alternative method. In
extending Larcher's proofs we have simplified his approach in a number
of ways. Most importantly we have avoided his use of the Dirichlet
theorem on primes in an arithmetic progression (for $\ZZ$).
\par
Strictly speaking, our definition of the cusp amplitude differs slightly
from the one given in [La]. The results in [La] are formulated for
$PSL_2(\ZZ)$, or equivalently, for congruence subgroups that contain
$-I_2$, and for those both definitions coincide.
\par
Cusp amplitudes were originally introduced for
subgroups of the modular group $SL_2(\ZZ)$ and play an important
role in the theory of modular forms and modular curves. Theorem A
holds (trivially) for every {\it normal} subgroup of $SL_2(D)$ since
then $\mathcal{A}(H)$ reduces to a single ideal. However it is known
that it does not hold in general for non-normal non-congruence
subgroups. Moreover we show that Theorem A does not hold in general
for the quasi-amplitudes of a congruence subgroup.
\par
The intersection of all the quasi-amplitudes (resp. cusp amplitudes)
of $H$ is $ql(H)$, its  quasi-level (resp. $l(H)$, its level).
As previously stated, for the classical case $D=\ZZ$ they clearly
coincide. Our second principal result shows, rather surprisingly,
that this is also true for ``many" congruence subgroups of $SL_2(D)$.
\par
We give a name to the precise condition as it will be used in several
places in the paper.
\\ \\
{\bf Condition L.}
Let $H$ be a subgroup of $SL_2(D)$ of non-zero level $l(H)=\qq$. We
say that $l(H)$ (or $H$) satisfies Condition L if both the following
two conditions hold:
\begin{itemize}
\item[(i)] $\qq+(2)=D$;
\item[(ii)] There are no prime ideals $\pp$ for which $|D/\pp|=3$ and
$\mathrm{ord}_{\pp}(\qq)=1$.
\end{itemize}
\noindent
Note that both conditions are satisfied if, for example, $\qq+(6)=D$,
so in particular always for $char(D)=p\geq 5$. Condition L is also
automatically satisfied if $D$ contains a field that is a nontrivial
extension of $\FF_3$.
\\ \\
{\bf Theorem B.} \it
Let $H$ be a congruence subgroup of $SL_2(D)$ whose level $l(H)=\qq$
satisfies Condition L. Then
$$ ql(H)=l(H).$$
In other words, then the quasi-level of $H$ is automatically an ideal.
\rm
\\ \\
There are examples of (normal) congruence subgroups which show that both
restrictions on $\qq$ from Condition L are necessary in Theorem B.
For non-congruence subgroups $H$ Theorem B breaks down completely,
i.e. $ql(H)$ and $l(H)$ can ``differ by as much as possible". More
precisely we prove that $SL_2(k[t])$, where $k$ is a field, contains
normal, non-congruence subgroups whose quasi-level has
$k$-codimension 1 in $k[t]$ and whose level can take any possible
value (including zero). (There are no proper subgroups of
$SL_2(k[t])$ whose quasi-level is $k[t]$.)
\par
Presumably one of the most interesting and in practice one of the most
useful results on congruence subgroups would be to bound the level in
terms of the index of the subgroup, i.e. to generalize the classical
result that the level of a congruence subgroup of $SL_2(\ZZ)$ can not
be bigger than its index. In Section 5 we discuss different
generalizations of this. The proofs use many of the previous results.
We highlight a version that holds for general $D$.
\\ \\
{\bf Theorem C.} \it
Let $D$ be any Dedekind domain and let $H$ be a congruence subgroup
of index $n$ in $SL_2(D)$ whose level $l(H)$ satisfies Condition L.
Then
$$|D/l(H)|\ \hbox{\it divides}\ n!,$$
and if $H$ is a normal subgroup then even
$$|D/l(H)|\ \hbox{\it divides}\ n.$$
\rm
\\
The structure of the paper is as follows:
\par
In Section 1 we clarify by means of theorems and counterexamples how
the properties `Every subgroup of finite index has non-zero level',
`Every subgroup of non-zero level is a congruence subgroup' and
`Every subgroup of finite index is a congruence subgroup' are related.
In Section 2 we present the quite lengthy proof of Theorem A.
In Section 3 we introduce quasi-amplitudes.
In Section 4 we define the quasi-level and prove Theorem B.
Using the previous results, we then can in Section 5 prove Theorem C
and similar relations between the level and the index of a congruence
subgroup.
\par
On the way we also prove several other interesting results.
Whenever appropriate we show by means of examples the necessity of
conditions in theorems and the limitations of certain definitions.
\par
Of course, every result that holds for congruence subgroups can also be
reformulated as a necessary condition that a subgroup must satisfy in
order to stand a chance of being a congruence subgroup. For the classical
modular group $SL_2(\ZZ)$ many such criteria are known, and some of them
are very useful in practice (compare Example 2.12).
In order to get an overview, we now summarize the criteria that our
results furnish. Several of them seem to be new or at least new in
such a general setting.
\\ \\
{\bf Theorem D.} \it
Let $D$ be any Dedekind domain and $\qq$ a non-zero ideal of $D$. Let $H$
be a subgroup of $G=SL_2(D)$ with level $l(H)=\qq$ and quasi-level $ql(H)$.
In order for $H$ to be a congruence subgroup, the following conditions are
necessary:
\begin{itemize}
\item[a)] (Corollary 1.4) \ $G(\qq)\subseteq H$.
\item[b)] (Corollary 2.8) \ The ideal-theoretic sum of all cusp amplitudes
of $H$ is again a cusp amplitude of $H$.
\item[c)] (Theorem 2.10) \ The intersection of all cusp amplitudes of
$H$ is again a cusp amplitude of $H$.
\item[d)] (Corollary 2.11) \ The level $\qq$ of $H$ is a cusp amplitude
of $H$.
\item[e)] (Lemma 4.2) \ For any $\alpha\in D$ that is invertible modulo
$\qq$ we have $\alpha^2 ql(H)\subseteq ql(H)$.
\end{itemize}
Now suppose moreover that the level $\qq$ of $H$ satisfies Condition L.
Then the following conditions are necessary for $H$ to be a congruence
subgroup:
\begin{itemize}
\item[f)] (Theorem 4.6) \ $ql(H)=l(H)$.
\item[g)] (Corollary 5.7) \ $|D/l(H)|\ \hbox{\it divides}\ (|G:H|)!$.
\end{itemize}
Finally suppose that in addition $H$ is a normal subgroup of $G$ (and
that $l(H)$ satisfies Condition L). Then the following conditions are
necessary for $H$ to be a congruence subgroup:
\begin{itemize}
\item[h)] (Corollary 4.7) \ Every quasi-amplitude of $H$ is actually
a cusp amplitude of $H$.
\item[i)] (Theorem 5.6) \ $|D/l(H)|\ \hbox{\it divides}\ |G:H|$.
\end{itemize}
\rm With the exception of a) in general none of the conditions is
sufficient for $H$ to be a congruence subgroup. Conditions e), f),
and h) hold automatically if $D=\ZZ$ and Remark 2.12 points out
non-congruence subgroups of $SL_2(\ZZ)$ that satisfy conditions b),
c), d), g) and i).
\par
On the other hand, none of the conditions is trivial in the sense of
automatically true for every $D$. Remark 2.12 provides examples of
non-congruence subgroups of $SL_2(\ZZ)$ that violate conditions
b), c) and d), and Theorem 4.12 allows the construction of normal
non-congruence subgroups that satisfy Condition L, but violate the
conditions e), f), g), h) and i).
\par
If $H$ is a congruence subgroup such that $l(H)$ does not satisfy
Condition L, then Examples 4.10 and 5.4 show that none of the conditions
f), g), h) and i) need hold, even if $H$ is normal. Compare also
Example 4.9 for conditions f) and h).
\par
Finally, the condition that $H$ is normal cannot simply be dropped
for parts h) and i). Examples 3.2 and 3.5 exhibit non-normal
congruence subgroups $H$ that satisfy Condition L, but neither h)
nor i).
\\

\subsection*{1. Congruence subgroups}

\noindent Let $R$ be a (commutative) ring and let $R^*$ be its
group of units. For each $r\in R,\alpha\in R^*$, we put
$$T(\alpha,r):=\left[\begin{array}{cc}\alpha&r\\
0&\alpha^{-1}\end{array}\right] \;\mathrm{and}\;
T(r):=T(1,r)= \left[\begin{array}{cc}1&r\\
0&1\end{array}\right].$$ In addition we define
$$S(r):=\left[\begin{array}{cc}1&0\\r&1\end{array}\right]\;
\mathrm{and}\; R(r):=\left[\begin{array}{cc}
1+r&r\\-r&1-r\end{array}\right].$$
\\
\noindent Let $E_2(R)$ be the subgroup of $SL_2(R)$ generated by all
$T(r)$ and $S(r)$ with $r \in R$. For each $R$-ideal $\qq$ we denote 
by $E_2(R,\qq)$ the {\it normal} subgroup of $E_2(R)$ generated by all
$T(q)$ with $q\in\qq$.  By definition we note that
$E_2(R,R)=E_2(R)$. Let $NE_2(R,\qq)$ be the {\it normal} subgroup of
$SL_2(R)$ generated by $E_2(R,\qq)$. We put $NE_2(R,R)=NE_2(R)$. We
also define
$$ SL_2(R,\qq)=\left\{ X \in SL_2(R):\:X\equiv I_2\;(\mod \qq)\right\}.$$
where $I_2$ denotes the $2$-dimensional unit matrix. Obviously
$$E_2(R,\qq)\leq NE_2(R,\qq)\trianglelefteq SL_2(R,\qq).$$
Occasionally we will also need the Borel group, consisting of the
upper triangular matrices,
$$B_2(R)=\{T(\alpha,r):\alpha \in R^*,r\in R\}.$$
If $H$ is a subgroup of $G$, we denote the core of $H$ in $G$, i.e. the
biggest {\it normal} subgroup of $G$ contained in $H$ by\\
$$N_H=\displaystyle{\bigcap_{g \in G}}H^g.$$
\\ \\
{\bf Definition.} The {\bf level} of a subgroup $H$ of
$SL_2(R)$, $l(H)$, is the largest ideal $\qq'$, say, for which
$NE_2(R,\qq')\leq H$. The level is well-defined since
$NE_2(R,\qq_1)\centerdot NE_2(R,\qq_2)=NE_2(R,\qq_1+\qq_2)$.
\par
Since $NE_2(R,\qq')$ is the normal subgroup of $SL_2(R)$ generated 
by all $T(r)$ with $r\in\qq'$, we can equivalently say that the level 
of $H$ is the largest ideal $\qq'$ such that the core $N_H$ of $H$ 
contains all translation matrices $T(r)$ with $r\in\qq'$.
\\ \\
{\bf Definition.} A subgroup $C$ of $SL_2(R)$ is called
a {\bf congruence subgroup} if
$$SL_2(R,\qq')\leq C$$
for some $\qq'\neq\{0\}$.
\\ \\
It is clear that if $R/\qq'$ is finite then $C$ is of finite index
in $SL_2(R)$.
\\ \\
\noindent Let $\qq,\;\qq'$  be ideals for which
$\qq\supseteq\qq'\supseteq\qq^2$. We will require some properties of the
quotient group $SL_2(R,\qq)/SL_2(R,\qq')$. We put $\rr=\qq/\qq'$
and we denote the image of any $q\in\qq$ in $\rr$ by
$\overline{q}$. Let $X\in SL_2(R,\qq)$. Then$$
X=\left[\begin{array}{cc} 1+x&y\\z&1+t\end{array}\right],$$ for
some $x,y,z,t\in\qq$. We define a
map$$\theta:SL_2(R,\qq)\longrightarrow \rr^3$$
 by $$\theta(X)=(\overline{x},\overline{y},\overline{z}).$$

 \noindent (Note that $x+t\equiv 0\;(\mod \qq').$)
\\ \\
{\bf Lemma 1.1.} \it With the above notation, the map
 $\theta$ induces the following isomorphism
 $$ SL_2(R,\qq)/SL_2(R,\qq')\cong (\rr^+)^3,$$
 where $\rr^+$ is the additive group of $\rr$.

 Moreover the quotient group $SL_2(R,\qq)/SL_2(R,\qq')$ is
 generated by the images of the elements $S(q),T(q),R(q)$, where
 $q\in\qq$.\\

\noindent {\bf Proof.} \rm See the proof of [MSt, Theorem 4.1].
\hfill $\Box$
\\ \\
To obtain more precise properties of congruence
 subgroups we require some restrictions on $R$. We recall that $R$
 is said to be an $SR_2$-{\it ring} if, for all $a,b\in R$ such that
 $aR+bR=R$, there exists $c\in R$ for which $(a+cb)\in R^*$.
Every semi-local ring, for example, is an
 $SR_2-$ring ([B, Theorem 3.5, p.239]).
\\ \\
{\bf Theorem 1.2.} \it Let $\qq,\qq'$ be $R$-ideals such
 that $R /\qq'$ is an $SR_2$-ring. Then
 $$ E_2(R,\qq)\centerdot SL_2(R,\qq')=SL_2(R,\qq+\qq').$$\\
 \noindent {\bf Proof.} \rm See the proof of [B, (9.3) Corollary,
 p. 267].
\hfill $\Box$
\\ \\
\noindent Since every proper quotient of a Dedekind domain is semi-local
the following is an immediate consequence of Theorem 1.2.
\\ \\
{\bf Corollary 1.3.} \it
Let $D$ be a Dedekind domain and let $\qq,\qq'$ be $D$-ideals.
 \begin{itemize}

 \item[(i)] If $\qq'\neq\{0\}$ then $$E_2(D,\qq)\centerdot SL_2(D,\qq')=SL_2(D,\qq+\qq').$$

 \item[(ii)] $$SL_2(D,\qq)\centerdot
 SL_2(D,\qq')=SL_2(D,\qq+\qq').$$
 \end{itemize}
 \noindent \rm The classical version of our next result (for the special case
 $D=\ZZ$) is due to Fricke. (See [W].)
\\ \\
{\bf Corollary 1.4.} \it
Let $D$ be a Dedekind domain and let $H$ be a subgroup of $SL_2(D)$ with
$l(H)=\qq\neq(0)$. Then $H$ is a congruence subgroup if and only if
$$H \geq SL_2(D,\qq).$$
\\
{\bf Remarks 1.5.} \rm
\begin{itemize}
\item[a)]
By virtue of Theorem 1.2, Corollaries 1.3 and 1.4 also apply to the case
where $R$ is a Noetherian domain of Krull dimension one.
For example all orders in algebraic number fields are of this type. One
``standard" example of a Noetherian domain of Krull dimension one which
is {\it not} Dedekind is the ring  $\ZZ[\sqrt{-3}]$. A Dedekind domain
$D$ is an {\it integrally closed} Noetherian domain of Krull dimension one.
\item[b)]
But Corollaries 1.3 and 1.4 do in general not hold for an integrally closed
Noetherian domain of Krull dimension bigger than one. Take for example
$R=\FF_q[x,y]$, the polynomial ring in two variables over a finite field.
Then
$$SL_2(R)/SL_2(R,yR)\cong SL_2(\FF_q[x]).$$
Let $f(x)\in\FF_q[x]$ be irreducible of degree at least $2$. By
Theorem 4.12 below or by [MSch2, Theorem 5.5] then $SL_2(\FF_q[x])$ has
a normal non-congruence subgroup $N$ of finite index and level $(f(x))$.
The inverse image of $N$ under reduction modulo $y$ is a normal subgroup
of finite index in $SL_2(R)$ of level $yR+f(x)R$ that contains
$SL_2(R,yR)$ but not $SL_2(R,yR+f(x)R)$.
\end{itemize}

\noindent
From now on $D$ will always be a Dedekind domain. We will sometimes
use the abbreviations
$$G=SL_2(D)$$
and
$$G(\qq)=SL_2(D,\qq)$$
for each $D$-ideal $\qq$.
\\ \\
{\bf Lemma 1.6.} \it
If $D$ has characteristic $p>0$ and every unit of $D$ has finite order,
then $D^*\cup\{0\}$ is a field $k$. Moreover, $k$ is the biggest algebraic
extension of $\FF_p$ contained in the field of fractions of $D$.
\\ \\
{\bf Proof.} \rm
Obviously $K$, the quotient field of $D$, contains the prime field
$\FF_p$. Let $k$ be the biggest algebraic extension of $\FF_p$
contained in $K$.
\par
Let $u$ be a unit of $D$ that has finite order. Then $u$ is a root
of unity and hence contained in an algebraic extension of $\FF_p$.
Thus $u\in k$. Conversely, every element of $k^*$ is a root of unity,
so it is in $D$ since $D$ is integrally closed, and obviously it is
a unit of $D$.
\hfill $\Box$
\\ \\
{\bf Lemma 1.7.} \it
If $D$ contains an infinite field, then every subgroup of finite index
in $SL_2(D)$ has level $D$.
\\ \\
{\bf Proof.} \rm
The following proof is inspired by the proof of [Se1, Proposition 1, p.491].
\par
If $H$ is a subgroup of finite index in $SL_2(D)$, then its core is
a normal subgroup of $SL_2(D)$ and still has finite index. So it
suffices to prove the lemma for a normal subgroup $N$ of finite
index in $SL_2(D)$. We define
$$ql(N):=\{r\in D\ :\ T(r)\in N\}.$$
This set will be investigated in more detail in Section 4. Here we only
need that $ql(N)$ is a subgroup of $(D,+)$ and that
$u^2 ql(N)\subseteq ql(N)$ for every $u\in D^*$. The last claim can be
seen by conjugating with the diagonal matrix with entries $u$ and
$u^{-1}$, as $N$ is normal.
\par
Let $k$ be the infinite field in $D$. Then obviously $k^*\subseteq D^*$.
So $ql(N)$ is stable under multiplication with $x^2$ and even with
$x^2 -y^2$ for all $x,y\in k$.
\par
If the characteristic of $k$ is different from $2$, then it is easy to
see that every element of $k$ is of the form $x^2 -y^2$ with $x,y\in k$.
Thus $ql(N)$ is a subspace of the $k$-vector space $D$. Since the quotient
space $D/ql(N)$ is finite because
$$[D:ql(N)]\le [SL_2(D):N]<\infty,$$
but $k$ is infinite, we must have $ql(N)=D$.
\par
If the characteristic is $2$, the squares in $k$ still form an infinite
field $k_2$, and we can apply the same proof using $k_2$-vector spaces.
\hfill $\Box$
\\ \\
{\bf Proposition 1.8.} \it
If $D$ has characteristic $0$ or if $D^*$ is infinite, then every
subgroup of finite index in $SL_2(D)$ has non-zero level.
\\ \\
{\bf Proof.} \rm
As in the previous proof we can assume that the subgroup $N$ is normal.
Again we show that $ql(N)$ contains a non-zero ideal. The first two cases
were already treated in the proof [Se1, Proposition 1, p.491], namely:
\par
If $char(D)=0$ and $[SL_2(D):N]=n$, then $T(nr)=(T(r))^n$ lies in $N$
for every $r\in D$, so $ql(N)$ contains the ideal $nD$.
\par
If $char(D)=p>0$ and $D^*$ contains a unit $u$ of infinite order, then
$D$ and $ql(N)$ are $\FF_p[u^2]$-modules. So the quotient group $D/ql(N)$
is also an $\FF_p[u^2]$-module, and since it is finite (whereas
$\FF_p[u^2]$ is isomorphic to a polynomial ring), it is annulated by
some nonzero element $a$ of $\FF_p[u^2]$. Thus $aD\subseteq ql(N)$.
\par
Finally, if $char(D)=p>0$ and there are no units of infinite order, then
Lemma 1.6 implies $D^*=k^*$ where $k$ is an infinite algebraic extension
of $\FF_p$. So in this case every subgroup of finite index has level $D$
by Lemma 1.7.
\hfill $\Box$
\\ \\
We recall from the introduction that if $D$ is of arithmetic type then
every congruence subgroup of $SL_2(D)$ is of finite index. It is known
that for some such $D$ (including the classical case $R=\ZZ$) the
converse does not hold. At this point we record a class of arithmetic
Dedekind domains for which the converse does hold.
\par
Let $K$ be a global field. In other words, $K$ is either an algebraic
number field, i.e. a finite extension of $\QQ$, or $K$ is an algebraic
function field of one variable with finite constant field, i.e., $K$ is
a finite extension of a rational function field $\FF_q(t)$.
\par
Let $S$ be a proper subset of all the places of $K$ which contains all
archimedean places. Let $\OOO_S$ be the set of all elements of $K$ that
are integral outside $S$. Then $\OOO_S$ is an arithmetic Dedekind domain
whose prime ideals correspond to the places of $K$ outside $S$.
\par
Note that we do not assume that $S$ is finite (as in the
standard definition of the {\it ring of $S$-integers} of $K$). If
$S$ is infinite, then $\OOO_S^*$ is not finitely generated. If $S$
contains all places of $K$ except one, then $\OOO_S$ is a discrete
valuation ring. The next theorem is a minor extension of a famous
result of Serre [Se1].
\\ \\
{\bf Theorem 1.9.}  \it
Let $K$ be a global field and $\OOO_S$ the ring of $S$-integers
where $S$ is an infinite set. Then every subgroup of finite index
in $SL_2(\OOO_S)$ has non-zero level, and every subgroup of non-zero
level is a congruence subgroup.
\rm
\\ \\
{\bf Proof.} The first claim follows from Proposition 1.8.
The second claim is equivalent to showing
$$SL_2(\OOO_S ,\qq)=NE_2(\OOO_S ,\qq)$$
for every non-zero ideal $\qq$ of $\OOO_S$. We use a ``local" argument
to reduce this to a suitable case where $S$ is finite. Take any
$g\in SL_2(\OOO_S,\qq)$. Then the entries of $g$, being elements of
$K$, lie in $\OOO_{S'}$ for some {\it finite} subset $S'$ of $S$. We
may assume that $|S'|>1$ and that $S'$ contains at least one
non-archimedean place. Let $\qq'$ be the ideal
$\qq\cap\OOO_{S'}$ of $\OOO_{S'}$. Then
$$SL_2(\OOO_{S'},\qq')=NE_2(\OOO_{S'},\qq')$$
by [Se1, Th\'eor\`eme 2 (b), p.498].
Since obviously $g\in SL_2(\OOO_{S'},\qq')$ and
$NE_2(\OOO_{S'},\qq')\subseteq NE_2(\OOO_S ,\qq)$, we have shown
$SL_2(\OOO_S ,\qq)\subseteq NE_2(\OOO_S ,\qq)$. The converse inclusion
always holds.
\hfill $\Box$
\\ \\
We present another class of Dedekind domains that have the same
congruence subgroup property.
\\ \\
{\bf Theorem 1.10.}  \it
If $D$ has only finitely many maximal ideals, then every subgroup
of finite index in $SL_2(D)$ has nonzero-level, and every subgroup
(finite index or not) of non-zero level is a congruence subgroup.
\rm
\\ \\
{\bf Proof.} \rm
Let $\aa$ be the product of the maximal ideals of $D$. Since the
ideals $\aa$, $\aa^2$, $\aa^3,\ldots$ are all different, $\aa$ must
be an infinite set. Now every element of the form $1+a$ with $a\in\aa$
obviously does not lie in any of the maximal ideals; so it must be
a unit. Thus $D^*$ is infinite. So if $H$ is a subgroup of finite
index in $SL_2(D)$, it has non-zero level by Proposition 1.8.
\par
Now let $H$ be a subgroup, not necessarily of finite index, and assume
that $H$ has non-zero level, say $l(H)=\qq$.
As $D$ is semi-local and hence an $SR_2$-ring by [B, Theorem 3.5, p.239],
we can apply Theorem 1.2 with $\qq'$ being the zero ideal, and obtain
$SL_2(D,\qq)\leq H$.
\hfill $\Box$
\\ \\
In the previous two results the reason for the congruence subgroup
property was that $NE_2(D,\qq)=SL_2(D,\qq)$ for every ideal $\qq$ of
$D$. We now construct examples with a different reason.
First a very general fact.
\\ \\
{\bf Lemma 1.11.}  \it
If $NE_2(D)\neq SL_2(D)$, then $NE_2(D,\qq)\neq SL_2(D,\qq)$ for
every non-zero ideal $\qq$ of $D$.
\rm
\\ \\
{\bf Proof.} \rm
Let $NE_2(D)\neq SL_2(D)$. Then $D$ has infinitely many maximal ideals
by Theorem 1.10. Asssume $NE_2(D,\qq)=SL_2(D,\qq)$. Choosing an ideal
$\qq'$ with $\qq +\qq'=D$ and using Corollary 1.3 we would obtain
$$NE_2(D)=NE_2(D,\qq+\qq')=NE_2(D,\qq)\centerdot NE_2(D,\qq')$$
$$=SL_2(D,\qq)\centerdot NE_2(D,\qq')=SL_2(D,\qq +\qq')=SL_2(D),$$
a contradiction.
\hfill $\Box$
\\ \\
{\bf Example 1.12.}
Let $k$ be an algebraically closed field of characteristic $0$.
Let $D=k[x,y]$ with $y^2 =x^3 +Ax+B$ such that the cubic polynomial
$x^3 +Ax+B\in k[x]$ has no multiple roots.
\par
Then $NE_2(D,\qq)\neq SL_2(D,\qq)$ for every non-zero ideal $\qq$ of
$D$. So $SL_2(D)$ has non-congruence subgroups of every level. But
every subgroup of finite index is a congruence subgroup for the trivial
reason that $SL_2(D)$ has no subgroups of finite index.
\par
To see all this, we use Takahashi's description [Ta] of the action of
$GL_2(D)$ on the Bruhat-Tits tree of $GL_2(k((t)))$ where $t=\frac{x}{y}$.
Note that since every element of $k^*$ is a square, there exists a
natural isomorphism $PSL_2(D)\cong PGL_2(D)$.
\par
Since the fundamental domain of this action is a connected graph
without loops [Ta], by the theory of groups acting on trees [Se2],
the group is generated by the stabilizers of the vertices. By [Ta,
Theorem 5] these stabilizers are built from subgroups isomorphic to
$k$, $k^*$ and $PSL_2(k)$. As the groups $k$ and $k^*$ are
infinitely divisible, they have no subgroups of finite index; and
since $k$ is infinite, $PSL_2(k)$ also has no finite index
subgroups. Hence a subgroup of finite index in $PGL_2(D)$ would have
to contain all stabilizers and so be equal to $PGL_2(D)$.
\par
Moreover, the fact that the fundamental domain has a vertex with trivial
stabilizer and the description of its neighbours [Ta] shows that
$PGL_2(D)$ is a free product
$$PGL_2(D)\cong\Star_{\ell\in\PP^1(k)}\Delta(\ell)$$
where $\Delta(\infty)$ is $GL_2(k)\star_{B_2(k)}B_2(D)$ modulo the center. 
This shows that modulo scalar matrices $SL_2(D)/NE_2(D)$ is isomorphic 
to the infinite group $\Star\limits_{\ell\in k}\Delta(\ell)$.
\\ \\
To appreciate this example we point out that even if a Dedekind domain
$D$ has no finite quotients, this does not imply that the group $SL_2(D)$
has no subgroups of finite index.
\\ \\
{\bf Example 1.13.}
Let
$$k=\bigcup_{n\in\NN}\FF_{2^{5^n}}\ \hbox{\rm and}\ D=k[x,y]\
\hbox{\rm with}\ y^2 +y=x^3.$$
Then $D^* =k^*$, and since every element of $k^*$ is a square we again
have $PSL_2(D)\cong PGL_2(D)$. By Takahashi's results [Ta] again we
have a free product
$$PSL_2(D)\cong\Star_{\ell\in\PP^1(k)}\Delta(\ell)$$
with $\Delta(\infty)$ equal to $SL_2(k)\star_{B_2(k)}B_2(D)$ modulo 
scalar matrices. However, since the equation $y^2 +y=1$ has no solution 
in $k$, by [Ta, Theorem 5] this time
$$\Delta(1)\cong L^*/k^*\ \hbox{\rm where}\ L=k(\omega)\ \hbox{\rm with}
\ \omega^2 +\omega=1.$$
So $L=\FF_4 k$ is the unique quadratic extension of $k$. Every element
of $k^*$ is a third power of an element of $k^*$. Moreover,
$\FF_8\nsubseteq L$, so $\FF_{64}\nsubseteq L$. Thus $L$ contains the
$3$-rd roots of unity, but not the $9$-th roots of unity. Hence the third
powers in $L^*$ form a subgroup of index $3$ that contains $k^*$. So
$L^*/k^*$ has a subgroup of index $3$. Because of the free product there
exists a surjective homomorphism from $PSL_2(D)$ to $\Delta(1)$ whose
kernel is the normal subgroup generated by all other $\Delta(\ell)$.
Thus $SL_2(D)$ contains a normal non-congruence subgroup of index $3$ and
level $D$.
\\ \\
We mention two types of Dedekind domains for which $SL_2(D)$ does not
have the congruence subgroup property.
\par
Every subgroup of finite index in $SL_2(\ZZ)$ has non-zero level, but
$SL_2(\ZZ)$ has finite index non-congruence subgroups.
\par
Let $\OOO_S$ be the ring of $S$-integers in a global function field
$K$ where $|S|=1$. Then $SL_2(\OOO_S)$ has uncountably many finite index
subgroups of level zero [MSch2, Corollary 3.6] and $SL_2(\OOO_S)$ has finite
index non-congruence subgroups of almost every level [MSch2, Theorem 5.5].
\\ \\
{\bf Remark 1.14.}
This leaves open the question whether there exist examples $SL_2(D)$ that
contain finite index subgroups of level zero but in which every finite
index subgroup of non-zero level is a congruence subgroup.
\par
By Lemmas 1.6 and 1.7 and Proposition 1.8 such a $D$ must be an
$\FF_q$-algebra with unit group $\FF_q^*$.
One might be tempted to think that the rings $\OOO_S$ in a global function
field $K$ with $|S|=1$ are the only instances of such rings. This is true
for finitely generated $\FF_q$-algebras $D$, because then the Krull
dimension of $D$ equals the transcendence degree of its quotient field.
\par
However, O.~Goldman [Go] has constructed Dedekind domains $D$ with
$D^*=\FF_q^*$ such that $D/\mm$ is a finite field for every maximal
ideal $\mm$ of $D$, and the quotient field $K$ of $D$ has
transcendence degree $n$ over $\FF_q$ where $n$ can be any given
natural number. We thank Peter V\'amos for pointing out this
reference to us.
\\

\subsection*{2. Cusp amplitudes}

\rm Throughout $D$ denotes a {\it Dedekind domain} with quotient
field $F (\neq D)$. For each $x \in D$ we denote the principal
$D$-ideal $xD$ by $(x)$. We write $(x,y)=z$,
if $(x)+(y)=(z)$. All the results of Section $1$ apply to $SL_2(D)$
which from now on we denote by $G$. \\
\noindent We recall the action of $G$ on $\hat{F}=\PP^1(F)$. Let
$g=\left[\begin{array} {cc}a&b\\c&d\end{array}\right]$
be an element of $G$. Then
$$g(z)=\left\{\begin{array} {cll}(az+b)/(cz+d)&,&z\in
F,\;(cz+d)\neq 0\\\infty&,&c\neq 0,\;z=-dc^{-1}\\ac^{-1}&,&c\neq 0,
\;z=\infty\\\infty&, &c=0,\;z=\infty
\end{array} \right. $$

\noindent In particular $G$ acts {\it transitively} on $\hat{F}$ if
and only if $D$ is a {\it principal ideal domain}. In general the
$G$-orbit in $\hat{F}$ containing $\infty$ is
$$O_{\infty}=G(\infty)= \{a/c:a,c
\in D,\;c\neq 0,\; (a,c)=1\}\cup\{\infty\}.$$\\
\noindent Let $H$ be a subgroup of $G$. For each $g\in G$, let
$H^g$ denote the conjugate subgroup $g^{-1}Hg$. For each
$z \in \hat{F}$, we denote the {\it stabilizer} of $z$ in $H$ by
$$H_z=\{h \in H: h(z)=z\}.$$
\noindent  It is clear that
$$ H_{g(z)}=((H^g)_z)^{g^{-1}},$$
for all $g\in G$ and $z\in \hat{F}$, and that
$$ G_{\infty}=B_2(D)=\{T(\alpha,r):\alpha \in D^*,r\in
D\}.$$ The following correspondence is obvious.
\\ \\
{\bf Definition.} Let $z_1,z_2 \in
\hat{F}$. We write $$z_1\equiv z_2\;(\mod H) \Longleftrightarrow
z_2=h(z_1),$$ for some $h\in H$.
\\ \\
{\bf Lemma 2.1.} \it Let $r,s \in G$. Then
$$ r(\infty) \equiv s(\infty)\; (\mod H) \Longleftrightarrow
HrG_{\infty}=HsG_{\infty}.$$

\noindent\rm  Lemma 2.1 provides the following bijection:
$$ H\backslash O_{\infty} \longleftrightarrow H\backslash G /G_{\infty};$$
between the orbits of the $H$-action on $O_{\infty}$ and the
$(H,G_{\infty})$ double cosets in $G$.
We refer to the elements of $H \backslash O_{\infty}$ as the $H$-{\bf cusps}.\\

\noindent It is clear that the double coset $HgG_{\infty}$ is a
union of a set $\mathcal{S}$ of (right) $H$-cosets and that
 $$ \mathcal{S} \longleftrightarrow G_{\infty}/G_{\infty}\cap H^g\longleftrightarrow G_z/H_z,$$
 where $z=g(\infty)$. From now on let $\{g_{\lambda}:\lambda\in\Lambda\}$ be a complete set
of representatives for the double coset space $H\backslash
G/G_{\infty}$. Then $\{z_{\lambda} = g_{\lambda}(\infty): \lambda
\in \Lambda\}$ is a complete set of representatives for the
$H$-cusps in $O_{\infty}$, by Lemma 2.1. Our next lemma is an
immediate consequence of the above.
\\ \\
{\bf Lemma 2.2.} \it Let $H$ have finite index in $G$.
Then, with the above notation,
$$|G:H|=\sum_{\lambda\in\Lambda}|G_{z_\lambda}:H_{z_\lambda}|=\sum_{\lambda\in\Lambda}|G_\infty:(H^{g_\lambda})_\infty|.$$
\\ \\
\noindent \rm It follows from the above that, if $g \in G$ and
$z=g(\infty)$, then every {\it unipotent} matrix in $G_z$ is of the
form $gT(x)g^{-1}$, for some $x \in D$. It is obvious that $\{ x \in
D:gT(x)g^{-1} \in H\}$ is a subgroup of (the additive group of) $D$.
Now let $g_i \in G, (i=1,2)$. Suppose now that
$Hg_1G_{\infty}=Hg_2G_{\infty}$. Then $g_1=hg_2s$, for some $h \in
H$ and $s \in G_{\infty}$. It is clear that
$$ g_1T(x)g_1^{-1} \in H \Longleftrightarrow g_2T(\mu x)g_2^{-1} \in H,$$
where $\mu \in D^*$ is completely determined by $s$. We are now able to make
the following definition.
\\ \\
{\bf Definition.} For each $g \in G$, let
$\cc(H,g)$ denote the {\it largest} $D$-ideal $\qq$ with the property that
$$ gT(q)g^{-1} \in H,$$
for all $q \in \qq$.
\\

\noindent It is clear from the above that $\cc(H,g)=\cc(H,g')$, whenever
$g(\infty) \equiv g'(\infty)\;(\mod H)$. We call $\cc(H,g)$ the
{\bf cusp amplitude} of the $H$-cusp containing $g(\infty)$.\\

\noindent It is clear that, for
each $g \in G$, there exists a unique $\lambda'\in\Lambda$ such that
$$ \cc(H,g) = \cc(H,g_{\lambda'}).$$
\noindent {\bf Definition.} We denote a {\it complete set of cusp
amplitudes for} $H$ by
$$\mathcal{A}(H)=\{\cc(H,g_{\lambda}):\lambda \in \Lambda\}.$$
\\
\noindent The following is an immediate consequence of the
above.
\\ \\
{\bf Lemma 2.3.} \it For each subgroup $H$ of $G$,
$$l(H)=\displaystyle{\bigcap_{g \in
G}}\cc(H,g)=\displaystyle{\bigcap_{\lambda \in
\Lambda}}\cc(H,g_{\lambda}).$$
\\
{\bf Remarks 2.4.}\rm \begin{itemize}
\item[(i)] If $N$ is a normal subgroup then $\cc(N,g_1)=\cc(N,g_2)$,
for all $g_1,g_2 \in G$. In this case $\mathcal{A}(N)$ involves a
single $D$-ideal.
\item[(ii)] If $l(H)$ is non-zero (for example when $H$ is a
congruence subgroup), then $\mathcal{A}(H)$ involves only {\it
finitely many} ideals.
\item[(iii)] To avoid any possible confusion we emphasize that if $D$ is
not a principal ideal domain the cusps (and the cusp amplitudes) of $H$
we consider are only a part of the $H$-orbits on $\PP^1(F)$, namely those
contained in the $G$-orbit of $\infty$.
\end{itemize}

\noindent It is easily verified that cusp amplitudes are invariant
under conjugation as given below.
\\ \\
{\bf Lemma 2.5.} \it Let $ k \in G$. Then
$$\cc(H^k,k^{-1}g)=\cc(H,g),$$ for all $g\in G$.\\

\noindent \rm We require a more detailed description of the
unipotent matrices in each $G_z$, where $z \in O_{\infty}$. If
$z=g(\infty)$ then $gT(-x)g^{-1} \in G_z$ is of the form
$$ U(a,b;x)= \left[\begin{array}{cc}
1+xab&-xa^2\\xb^2&1-xab\end{array}\right],$$ where $(a,b)=1$. In
this case $z=a/b$, when $b\neq 0$, and $z=\infty$, when $b=0$.
Note that $U(1,0;x)=T(x)$ and $U(0,1;x)=S(x).$\\

\noindent Before our principal results we record this well-known
useful property of Dedekind domains, which follows from the Chinese
Remainder Theorem (CRT).
\\ \\
{\bf Lemma 2.6.} \it Let $\pp_1,\cdots,\pp_t$ be distinct
prime $D$-ideals and let $\alpha_1,\cdots,\alpha_t$ be non-negative
integers, where $t\geq 1$. Then there exists $d \in D$ such that

 $$d \in \pp_i^{\alpha_i}\backslash \pp_i^{\alpha_i+1},$$
 where $1\leq i \leq t.$\\

\noindent \rm For the first principal result our approach is more direct
than that of Larcher [La]. In particular we avoid his use of the
Dirichlet theorem on primes in an arithmetic progression
(for $\ZZ$.)
\\ \\
{\bf Theorem 2.7.} \it Let $H$ be a congruence subgroup of
$G$ and let $\qq_i$ be a non-zero $D$-ideal contained in
$\cc(H,g_i)$, where $i=1,2$. Then there exists $g_0 \in G$ such that
$$ \qq_1+\qq_2 \subseteq \cc(H,g_0).$$
\\
\noindent {\bf Proof.} \rm Now $H \geq G(\qq)$, for some $\qq\neq
\{0\}$. It is clear that, if $X,Y \in G$ and $X \equiv Y
(\mod \qq)$, then $X \in H$ if and only if $Y \in H$.\\
\noindent By Lemma 2.5 it is sufficient to prove the theorem for the
case $g_1=I_2$. We may also assume that $\qq_1+\qq_2 \neq \qq_i$ and
(by Corollary 1.4
and Lemma 2.3) that $\qq_i\supseteq\qq$, where $i=1,2$.\\
\noindent Now $U(1,0;q)=T(q) \in H$, for all $q \in \qq_1$ and
$U(c,d;q) \in H$, for all $q \in \qq_2$, for some $c,d \in D$, with
$(c,d)=1$. Let $\qq_1+\qq_2=\qq_0$. For our purposes it is
sufficient to find $a,b$, with $(a,b)=1$ and ideals
$\overline{\qq}_i$, contained in $\qq_i$, where $i=1,2$, such that
$\overline{\qq}_1+\overline{\qq}_2=\qq_0$, with the following
properties. For all $\lambda \in \overline{\qq}_1$ and all $\mu \in
\overline{\qq}_2$, there exist $\lambda_0 \in \qq_1$ and $\mu_0 \in
\qq_2$ such that:
\begin{itemize}
\item[(a)]$U(a,b;\lambda) \equiv T(\lambda_{0})\; (\mod \qq)$;
\item[(b)]$U(a,b;\mu) \equiv U(c,d;\mu_{0})\; (\mod
\qq).$
\end{itemize}
\noindent We may then take $g_0$ to be any matrix of the form
$$\left[\begin{array} {cc}a&*\\b&*\end{array}\right]\in G.$$
\noindent  By definition $\qq_i=\qq_i'\qq_0$, where $i=1,2$. Then
$\qq_1'+\qq_2'=D$. Let $\qq_i''$ be the ``smallest" divisor of
$\qq_0^{-1}\qq$ all of whose prime divisors also divide $\qq_i'$,
where $i=1,2$.
(Clearly $\qq_i'' \subseteq\qq_i'$.) Then $\qq_1''+\qq_2''=D$ and$$
\qq=\qq_1''\qq_2''\qq_0\rr,$$say, where $\rr+\qq_1''\qq_2''=D$.\\
\\
\noindent Since $(a,b)=1$ the congruence (a) is equivalent to
$$ b\lambda\equiv \lambda_0+a^2\lambda
\equiv 0\;(\mod \qq).\;\;\;\;\;\;\;\;(*)$$ \noindent Since $(c,d)=1$
the congruence (b) is equivalent to: $$ab\mu\equiv cd\mu x\;(\mod
\qq)\; ;\;a^2\mu\equiv c^2\mu x\:(\mod \qq)\;;\;b^2\mu\equiv d^2\mu
x\:(\mod \qq)\; ;\;(**)$$for some $x \in D$ such that $\mu_0\equiv
\mu x\;(\mod \qq)$. A consequence of $(**)$ is that, for all $\mu \in
\overline{\qq}_2$, $$ (ad-bc)\mu \equiv 0\; (\mod \qq).$$ \noindent
We begin by finding $a,b\in D$ for which \begin{itemize}
\item [(i)]$ ad \equiv bc\;(\mod \qq_1'')$,
\item[(ii)] $(a,b)=1$,
\item[(iii)] $b \in \qq_2''\rr$. \end{itemize}
\noindent By Lemma 2.6 we may choose $d_0 \in (c)+\qq_1''$ and $b
\in \qq_2''\rr$ such that \begin{itemize} \item[(i)]
$(dd_0)+\qq_1''=(bc)+\qq_1'',$ \item[(ii)]$(d_0)+(b)=D.$
\end{itemize}
\noindent By applying the CRT to the factors in the prime
decomposition of $\qq_1''$ we can find $y \in D$ such that
\begin{itemize}
\item[(i)] $ydd_0 \equiv bc\;(\mod \qq_1''),$
\item[(ii)] $(y)+\qq_1''=D.$
\end{itemize}
\noindent (If $\pp^{\alpha}||\qq_1''$, where $\pp$ is prime with
$\alpha >0$, and $dd_0\equiv bc \equiv 0\;(\mod \qq_1'')$, take
$y=1$.) Now for each prime divisor $\pp$ of $(b)$ we again use the
CRT to find $z \in D$ for which \begin{itemize}\item[(i)] $z \equiv
y\;(\mod \pp^{\alpha})$, where $\pp^{\alpha}||\qq_1''$, \item[(ii)]
$z \equiv 1\;(\mod \pp)$, when $\pp\nmid \qq_1''$. \end{itemize}
\noindent The elements $a=zd_0$ and $b$ (as above) satisfy the
requirements. \\
\noindent We can use the CRT to find $ x \in D$ such that
$$ ab \equiv cdx\:(\mod \qq_1'')\; ,\; a^2\equiv c^2x\;(\mod \qq_1'')\; ,\; b^2\equiv d^2x\;(\mod \qq_1'').$$
\noindent (Suppose that $\pp^{\alpha}||\qq_1''$. If $c \notin \pp$,
take $x'\equiv c^{-2}a^2\;(\mod \pp^{\alpha})$. If $d \notin \pp$,
take $x'\equiv d^{-2}b^2\;(\mod \pp^{\alpha})$.) Then, for all $\mu
\in
\overline{\qq}_2=\qq_2''\qq_0\rr$, the congruences $(**)$ are satisfied.\\
\noindent In addition, for all $\lambda \in
\overline{\qq}_1=\qq_1''\qq_0$, the congruences $(*)$ are satisfied
(with $\lambda_0=-a^2\lambda$). This completes the proof.
\hfill $\Box$
\\ \\
{\bf Corollary 2.8.} \it Let $H$ be a congruence subgroup
of $G$ and let
$$ \cc_{\mathrm{max}}= \displaystyle{\sum_{\qq \in
\mathcal{A}(H)}\qq}=\displaystyle{\sum_{g \in G}\cc(H,g)}.$$ Then
$$\cc_{\mathrm{max}} \in \mathcal{A}(H).$$
\\
{\bf Proof.} \rm As previously noted $\mathcal{A}(H)$
involves only finitely many distinct ideals $\qq_1\cdots\qq_t$, say.
By repeated applications of Theorem 2.7 it follows that there exists
$g' \in G$ such that
$$\qq_1+\cdots+\qq_t \subseteq \cc(H,g').$$
But $\cc(H,g')\in\mathcal{A}(H)$ and so
$\cc(H,g')=\qq_1+\cdots+\qq_t=\cc_{\mathrm{max}}$.
\hfill $\Box$
\\ \\
We now come to our second principal result. Here our
approach is similar to that of Larcher.
\\ \\
{\bf Lemma 2.9.} \it Let $H$ be a congruence subgroup of
$G$ of level $\qq$ and let $\qq_{\infty}= \cc(H,I_2)$ and
$\qq_0=\cc(H,g_0)$, where $g_0(\infty)=0$. Then
$$ \qq\supseteq\qq_0\qq_{\infty}.$$
\\
{\bf Proof.} \rm Suppose to the contrary that
$\qq_0\qq_{\infty}+\qq\neq \qq$. We will prove that there exists an
ideal $\qq'$ such that
$$ G(\qq')\leq H\;\mathrm{and}\;\qq'
\supsetneq \qq,$$ which contradicts Corollary 1.4.

\noindent Recall that $\qq_{\infty}\supseteq\qq$ and
$\qq_0\supseteq\qq$, by Lemma 2.3. We assume for now that
$$ \qq_{\infty}\qq_0+\qq=\qq',$$
where $\pp\qq'=\qq$, for some prime ideal $\pp$. Then $T(x),S(x) \in
H$, for all $x \in \qq'$. There are two
possibilities.
\\ \\
\noindent {\bf Case 1: $\pp+\qq'=D.$}
\\ \\
\noindent By Corollary 1.3
$$G(\qq')/G(\qq) \cong G/G(\pp)\cong SL_2(D/\pp).$$
Now $D/\pp$ is a field and so $SL_2(D/\pp)$ is
generated by elementary matrices. It follows that
$$ \langle T(x),S(x)\;(x \in \qq'),\;G(\qq)\rangle=G(\qq'),$$
and hence that $G(\qq') \leq H$.
\\ \\
\noindent {\bf Case 2: $\pp+\qq'=\pp.$}
\\ \\
\noindent Note that here Lemma 1.1 applies to the quotient group
$G(\qq')/G(\qq)$, since $(\qq')^2\subseteq\qq$. For all $x \in
\qq_{\infty}$ and all $y \in \qq_0$, the element
$$ V=T(x)S(y)T(-x)S(-y) \in H\cap
SL_2(D,\qq').$$ Now
$$ V \equiv \left[ \begin{array}{cc}1+xy&*\\
*&1-xy\end{array}\right]\;(\mod \qq).$$
\noindent It follows that $H\cap G(\qq')$ contains elements
$$W \equiv \left[ \begin{array}{cc}1+z & * \\ * &
1-z\end{array} \right]\;(\mod \qq),$$ for all $z \in \qq'$ and also
$S(z)$ and $T(z)$ (as mentioned before Case 1). This implies that
$G(\qq')\leq H$, by Lemma 1.1.
\\ \\
\noindent If $\qq_{\infty}\qq_0+\qq\neq \qq'$ we replace $\qq_0$,
say, with $\qq^*$, divisible by $\qq_0$, where
$\qq_{\infty}\qq^*+\qq=\qq'$, and repeat the above argument.
\hfill $\Box$\\ \\

\noindent An alternative proof of our second principal result (for
the case $D=\ZZ$) can be found in [St].
\\ \\
{\bf Theorem 2.10.} \it Let $H$ be a congruence subgroup of
$G$ and let
$$ \cc_{\mathrm{min}}= \displaystyle{\bigcap_{\qq \in
\mathcal{A}(H)}\qq}=\displaystyle{\bigcap_{g \in G}\cc(H,g)}.$$
\noindent Then
$$\cc_{\mathrm{min}} \in \mathcal{A}(H).$$

\noindent {\bf Proof.} \rm  Recall from Lemma 2.3 that
$\cc(H,g)\supseteq\qq$, for all $g \in G$, where $\qq=l(H)$. We will
assume that
$$ \cc(H,g) \neq \qq,$$
for all $ g \in G$, and obtain a contradiction.\\ \\
\noindent By Lemma 2.5 it is sufficient to prove the theorem for the
case where $$\cc_{\mathrm{max}}=\cc(H,I_2).$$ (See Corollary 2.8.)
Let $\cc(H,I_2)=\qq_{\infty}$ and let $\cc(H,g_0)= \qq_0$, where
$g_0 \in G$ is any element for which $g_0(\infty)=0$. ( For example,
$g_0=T(-1)S(1)$.) By Corollary 2.8 it follows that
$\qq_{\infty}\supseteq\qq_0$. Now
$\qq=\qq_{\infty}\qq_{\infty}'=\qq_0\qq_0'$, say. Let
$\rr=\qq_{\infty}+\qq_{\infty}'$. Then $\qq_0'\supseteq \rr$, by
Lemma 2.9.\\

\noindent By Lemma 2.6 we can choose $z \in D$ for prime ideals
$\pp$ such that
\begin{itemize}
\item[(i)] $z \in \pp$, when $\pp|\rr,\;\pp\nmid \qq_0'$,
\item[(ii)] $z \notin \pp$, when $\pp|\qq_0'$.
\end{itemize}(Possibly $z=1$.) Consider the $H$-cusp,
$\cc(H,g_z)$, where $g_z \in G$ is given by
$$g_z=\left[\begin{array}{cc}z&*\\ 1&*\end{array}\right].$$
Then $g_z(\infty)=z$ and $\cc(H,g_z)=\qq_z$ is determined by the
matrices
$$U(z,1;x)=Q(x),$$
where $x \in D$. Let $\qq=\qq_z\qq_z'$. If $t=T(z)$, then, by Lemma
2.5,
$$ \cc(H,I_2)=\cc(H^t,t^{-1})\;\mathrm{and}\;\cc(H,g_z)=\cc(H^t,t^{-1}g_z).$$
Now $t^{-1}(\infty)=\infty$ and $t^{-1}g_z(\infty)=0$. Applying
Lemma 2.9 therefore to $H^t$, it follows that $\qq_z'\supseteq \rr$.
(By definition $l(H)=l(H^t)=\qq$.) Note that by our initial
hypothesis $\qq_0',\qq_z' \neq D$. There are two possibilities.
\\ \\
\noindent {\bf Case 1: $\qq_0'+\qq_z'\neq D.$}
\\ \\
\noindent Choose a prime divisor $\pp$ of $\qq_0'+\qq_z'$. Let
$\qq=\pp\qq'$. Then $(\qq')^2\subseteq \qq$, since $\pp$ divides
$\qq_{\infty}$ and $\qq_{\infty}'$. It is clear that $\qq'$ is
divisible by $\qq_{\infty},\qq_0,\qq_z$. It follows that
$T(x),S(x),Q(x) \in H $, for all $ x \in \qq'$.\\ \\
\noindent Now by definition $\pp +(z)=D$ and so $\qq+(z)\qq'=\qq'$.
Hence, for all $q \in \qq'$, there exist $q' \in \qq'$ such that
$zq' \equiv q\;(\mod \qq)$. As in the proof of Lemma 2.9 (Case 2) it
follows from Lemma 1.3 that $G(\qq') \leq H$, which contradicts
Corollary 1.4.
\\ \\
\noindent {\bf Case 2: $\qq_0'+\qq_z'=D.$}
\\ \\
\noindent Choose a prime $\pp$ dividing $\qq_z'$. Since
$\qq_z' \supseteq\rr$ and $\pp+\qq_0'= D$, it follows that
$z \in \pp$. As before let $\qq=\pp\qq'$. Then $\qq'\subseteq\qq_z$,
so $Q(x) \in H$ and
$$ Q(x) \equiv S(x)\;(\mod \qq),$$ for all $x  \in \qq'$. It follows that $S(x) \in H$, for all $x \in
\qq'+\qq_0$. But $\qq'+\qq_0 \neq \qq_0$, since $\qq_0 \nsupseteq
\qq'$, which contradicts the maximality of
$\cc(H,g_0)$. This completes the proof.
\hfill $\Box$
\\ \\
From Theorem 2.10 we obtain in particular the following necessary
condition for a subgroup to be a congruence subgroup.
\\ \\
{\bf Corollary 2.11.} \it If $H$ is a congruence subgroup of $G$,
then there exists some $g\in G$ with $l(H)=\cc(H,g)$. \rm
\\ \\
\noindent All the results in this section are, of course, motivated
by the classical modular group. It is appropriate therefore that we
include a detailed account of how they apply to this important
special case.
\\ \\
{\bf Example 2.12.} Cusp amplitudes were first introduced
for finite index subgroups of the modular group $SL_2(\ZZ)$
(motivated by the theory of modular forms). In this case the cusp
$\infty$ is $\PP^1(\QQ)$, since $\ZZ$ is a PID. Let $H$ be a
subgroup  of {\it finite index} in $G$ containing $\pm I_2$. For
each $g \in G$, the $\ZZ$-ideal $\cc(H,g)$ is non-zero, generated by
a (unique) positive integer $n(H,g)$, say. Since $\ZZ^*=\{\pm 1\}$
it follows that in this case
$$|G_{\infty}:G_{\infty}\cap H^g|=|G_z:H_z|=n(H,g),$$
where $z=g(\infty)$. The double coset space $H \backslash
G/G_{\infty}$ is finite, represented by elements $g_1,\cdots,g_t$,
say, of $G$. Let $n_i = n(H,g_i)$, where $1\leq i\leq t$. We assume
that $n_1 \leq \cdots \leq n_t$. The sequence $(n_1,\cdots,n_t)$ is
called the {\it cusp-split} of $H$ and satisfies the important {\it
cusp-split equation}
$$ n_1+\cdots+n_t=|SL_2(\ZZ):H|.$$

\noindent If $n_0(\geq 1)$ is a generator of the ($\ZZ$-)ideal $l(H)$,
then $n_0=\mathrm{lcm}\{n_1,\cdots,n_t\}$, by Lemma 2.3.
\\ \\
\noindent Suppose now that $H$ is a congruence subgroup of index
$\mu$ in $SL_2(\ZZ)$. Then, by Corollary 2.8 and Theorem 2.10,
$$ n_1=\mathrm{gcd}\{n_1,\cdots,n_t\}\ \ \mathrm{and}\ \ n_t=n_0.$$
\noindent There are two immediate consequences, namely,
$$\mu \geq n_0,\qquad\qquad (*)$$
\noindent and
$$\mu\quad\mathrm{divides}\quad|SL_2(\ZZ):SL_2(\ZZ,n_0\ZZ)|.\qquad\qquad (**)$$
\noindent (See Corollary 1.4.) The inequality $(*)$ can be used [St] to
obtain an upper bound for the number of congruence subgroups of
$SL_2(\ZZ)$ of bounded index.
\\ \\
Among the index $7$ subgroups of $SL_2(\ZZ)$ are those
with cusp-splits $(3,4)$ and $(2,5)$. (See [St].) They are
non-congruence by Corollary 2.8 or Theorem 2.10 (or inequality $(*)$).
On the other hand these results provide only {\it necessary}
conditions. For example, it is known [AS] that there exist subgroups
with cusp-splits $(1,6)$ and $(1,1,7)$ which, despite being
consistent with Corollary 2.8 and Theorem 2.10, are non-congruence
since neither satisfies $(**)$. However non-congruence subgroups do
exist for which Corollary 2.8, Theorem 2.10 and $(**)$ all hold. It is
known [AS] there exists a non-congruence subgroup with cusp-split
$(8)$.
\\ \\
As previously emphasised  we wish to demonstrate
explicitly that our principal results, Corollary 2.8 and Theorem
2.10, do extend to $D$ which are not of arithmetic type.
\\ \\
{\bf Example 2.13.}  Let $D=k[t]$, the polynomial ring
over a field $k$. Then $D$ is of arithmetic type if and only if $k$
is {\it finite}. Let $\pp$ be
a prime $k[t]$-ideal.\\ \\
\noindent We define
$$ G_0(\pp)=\langle T(\alpha,\beta)\;(\alpha \in k^*,\beta\in
k),\;SL_2(k[t],\pp)\rangle.$$\\
\noindent Now $T(x) \in G_0(\pp)$ if and only if $x \in k+\pp$. Let
$g \in G$. Then $z=g(\infty)=a/b$, where $a,b\in k[t]$, with
$(a,b)=1$. If $U(a,b;x)\in G_0(\pp)$, then $xb^2\equiv 0\;(\mod
\pp)$. There are
two possibilities.\\ \\
 {\bf (i): $b \notin \pp$.}\\ \\
 In this case $x \in \pp$ and hence $\cc(G_0(\pp),g)=\pp$.\\ \\
 \noindent {\bf (ii): $b \in \pp$.}\\ \\
 \noindent In this case $a \notin \pp$ and $U(a,b;x) \equiv T(-xa^2)\;(\mod
 \pp)$. Thus $U(a,b;x)\in G_0(\pp)$ if and only if $ x\in c^2k+\pp$,
 where $ac\equiv 1\;(\mod \pp)$. We are reduced to two possible
 outcomes.\\  \\
 If $\dim_k(k[t]/\pp)>1$, then
 $$\cc(G_0(\pp),g)=\pp,$$for all $g \in G$. We note that
 $G_0(\pp)$ is an example of a {\it non-normal} subgroup all of whose
 cusp-amplitudes are equal.\\ \\
 \noindent Suppose then that $\dim_k(k[t]/\pp)=1$ (equivalently,
$k[t]/\pp \cong k$). (This always
 happens when $k$ is algebraically closed.) It is easily verified
 that there are precisely two $(G_0(\pp),S_{\infty})$ double cosets
 in $G=SL_2(k[t])$ and that $\{T(1),S(1)\}$ is a complete set of
 representatives for the $G_0(\pp)$-cusps, $G_0(\pp)\backslash
 \widehat{k(t)})$. It follows that
 $$ \mathcal{A}(G_0(\pp))=\{k[t],\pp\}.$$

\subsection*{3. Quasi-amplitudes}

\noindent {\bf Definition.} As usual let $H$ be a subgroup of $G$.
We define the {\bf quasi-amplitude}
$$\bb(H,g)=\{b\in D: gT(b)g^{-1}\in H\}=\{b\in D:T(b)\in H^g\}.$$
\\
Clearly the cusp amplitude $\cc(H,g)$ is the biggest $D$-ideal
contained in $\bb(H,g)$.
\par
\noindent Note that, in contrast to $\cc(H,g)$, the quasi-amplitude
$\bb(H,g)$ does not just depend on the
$H$-cusp containing $g(\infty)$. More precisely it is easily verified that\\
$$\mathrm{if}\; g({\infty}) \equiv g'({\infty})\;(\mod H),\;
\mathrm{then}\; \bb(H,g)=u^2\bb(H,g'),$$ for some $u \in D^*$. It is
clear that if $H$ is {\it normal} in $G$, then
$\bb(H,g_1)=\bb(H,g_2)$, for all $g_1,g_2 \in G$.
\par
\noindent We put
$$U=\{T(r): r\in D\}.$$
\noindent It is clear that, for all $g \in G$, the subgroup
$(H^g)_{\infty}$ normalizes $U$, so that $U\cdot(H^g)_{\infty}$ is a
subgroup of $G$. We then have the following one-one correspondences
$$ D/\bb(H,g)\leftrightarrow U/U\cap(H^g)_{\infty} \leftrightarrow U\cdot (H^g)_{\infty}/(H^g)_{\infty}.$$
\noindent {\bf Definition.} If $H$ is of finite index in $G$ we
define$$m(H,g)=|G_\infty :U\cdot(H^g)_\infty|.$$ \noindent Our next
result is a generalization of the cusp split formula for subgroups
of $SL_2(\ZZ)$.
\\ \\
{\bf Theorem 3.1.} \it Let $H$ be a subgroup of finite
index in $G$. If $\{g_\lambda(\infty): \lambda\in\Lambda\}$ is a
complete set of representatives for the $H$-cusps in $O_\infty$,
then
$$|G:H|=\sum_{\lambda\in\Lambda}m(H,g_\lambda)|D:\bb(H,g_\lambda)|.$$
\noindent {\bf Proof.} \rm  We note that, for all $g \in G$,
\\$$|G_\infty :(H^g)_\infty|=|G_\infty :U\cdot (H^g)_\infty|\cdot
|U\cdot (H^g)_\infty:(H^g)_\infty|=m(H,g)|D:\bb(H,g)|.$$ \noindent
The proof follows from Lemma 2.2.
\hfill $\Box$
\\ \\
{\bf Example 3.2.}
Let $\CCC=\FF_9[t]$. Then
$$SL_2(\CCC)/\pm SL_2(\CCC,(t))\cong PSL_2(\FF_9)\cong A_6.$$
So there exists a (non-normal) congruence subgroup $H$ of
$SL_2(\FF_9[t])$ of index $6$ and level $(t)$.
\par
Obviously $(t)\subseteq\bb(H,I_2)\subseteq\FF_9[t]$. Actually,
both inclusions are proper: $\bb(H,I_2)$ cannot equal $\FF_9[t]$
since $H/\pm SL_2(\CCC,(t))\cong A_5$ does not contain a subgroup
of order $9$; on the other hand
$$|\CCC:\bb(H,I_2)|\le |SL_2(\CCC):H|=6,$$
so $\bb(H,I_2)$ cannot equal $(t)$. In particular, this shows that
$\bb(H,I_2)$ is not an $\FF_9$-subspace of $\FF_9[t]$, although,
being an additive group, it is clearly an $\FF_3$-subspace.
All in all we see that
$$\bb(H,I_2)=\FF_3\beta+(t)$$
where $\beta$ is an element from $\FF_9^*$.
\par
Let $\zeta\in\FF_9$ be a primitive $8$-th root of unity. Then
$H$ cannot contain any matrices of the form
$$\left[\begin{array}{cc}
\zeta & *\\
0     & \zeta^{-1}\\
\end{array}
\right]$$
since $H/\pm SL_2(\CCC,(t))$ has no elements of order $4$. Hence
$|G_{\infty}:H_{\infty}|$ is divisible by $2$. Together with
$|\CCC:\bb(H,I_2)|=3$ we see that $|G_{\infty}:H_{\infty}|=6$
and hence that $\infty$ is the only cusp of $H$.
\par
Here we see the problem discussed earlier. The matrix
$$\gamma=\left[\begin{array}{cc}
\zeta & 0\\
0     & \zeta^{-1}\\
\end{array}
\right]$$
represents the same cusp as $I_2$, namely $\infty$, but
$$\bb(H,\gamma)=\zeta^2\bb(H,I_2)=\FF_3 i\beta+(t)\neq\bb(H,I_2)$$
where $i$ is a primitive $4$-th root of unity in $\FF_9$.
\\ \\
{\bf Proposition 3.3.} \it
\begin{itemize}
\item[a)]
If $q>3$, then $SL_2(\FF_q[t])$ has no proper, normal subgroup
of finite index with only one cusp.
\item[b)]
If $q\le 3$, then for every positive integer $e$ there are uncountably
many normal subgroups of index $q^e$ in $SL_2(\FF_q[t])$ that have
only one cusp.
\end{itemize}
\rm
\noindent
{\bf Proof.} \rm
a) If $N$ is a normal subgroup of finite index in $G$ with only one
cusp, then Theorem 3.1 shows that the Borel group $G_\infty$ contains
a system of coset representatives for $N$ in $G$. Thus $G/N$ is a
semidirect product of an elementary abelian $p$-group where
$p=char(\FF_q)$ and a cyclic group of order $m$ where $m|(q-1)$.
If $m>1$, then $G/N$ contains a normal subgroup of index $m$, in
contradiction to the fact that the minimal index of a normal subgroup
in $G$ is $|PSL_2(\FF_q)|$. Compare [MSch1, Theorem 6.2]. If $m=1$,
then $G/N$ contains a normal subgroup of index $p$, again contradicting
the same minimal index.
\par
b) Recall Nagao's Theorem
$$SL_2(\FF_q[t])=SL_2(\FF_q)\star_{B_2(\FF_q)}B_2(\FF_q[t])$$
Compare for example [Se2, exercise 2, p.88]. We define
$$V=\{T(r)\ :\ r\in (t)\}.$$
Since for $q\in\{2,3\}$ all diagonal matrices are in the center of
$SL_2(\FF_q[t])$, this amalgamated product shows that the identity on $V$
extends to a surjective group homomorphism $\phi: SL_2(\FF_q[t])\to V$
whose kernel is the normal subgroup generated by $SL_2(\FF_q)$. Now for
each of the uncountably many $\FF_q$-subspaces $W$ of codimension $q^e$
in $V$ the inverse image $\phi^{-1}(W)$ is a normal subgroup of
$SL_2(\FF_q[t])$ of index $q^e$, and Theorem 3.1 shows that it has only
one cusp.
\hfill $\Box$
\\ \\
\noindent The following lemma will be required later on.
\\ \\
{\bf Lemma 3.4.} \it Let $H$ be a congruence subgroup of
$G$ of level $\qq$. Assume that $x\in\bb(H,I_2)$. Then for every
$\alpha\in D$ that is invertible modulo $\qq$ there exists an
element $g\in G$ with $\alpha^2 x\in\bb(H,g)$. \rm
\\ \\
{\bf Proof.} \rm Since $\alpha$ is invertible mod $\qq$,
there exists $\delta\in D$ such that $\alpha\delta-1=:\gamma\in\qq$.
Let $g^{-1}={\alpha\ 1\choose\gamma\ \delta}\in G$. Then $g^{-1}Hg$
contains
$$g^{-1}T(x)g=\left[\begin{array}{cc}
1-\alpha\gamma x & \alpha^2 x\\
-\gamma^2 x        & 1+\alpha\gamma x\\
\end{array}\right].$$
Moreover,
$$g^{-1}T(x)gT(-\alpha^2 x)=\left[\begin{array}{cc}
1-\alpha\gamma x & \alpha^3\gamma x\\
-\gamma^2 x        & 1+\alpha^2\gamma^2 x^2+\alpha\gamma x\\
\end{array}\right]$$
is an element of $G(\qq)$, which is a subgroup of $g^{-1}Hg$. Hence
$T(\alpha^2 x)\in g^{-1}Hg$.
\par{}\quad
\hfill $\Box$
\\ \\
\noindent We conclude this section with two examples which show that
in general neither Corollary 2.8 nor Theorem 2.10 hold for the
quasi-amplitudes of a congruence subgroup. If $H$ is any subgroup of
$G$, we denote the subgroup $\langle H,-I_2\rangle$ by $\pm H$.
\\ \\
{\bf Example 3.5.}
Here we use $G$ to denote the {\it Bianchi group} $SL_2(\OOO_d)$,
where $\OOO_d$ is the ring of integers in the imaginary quadratic
number field $\QQ(\sqrt{-d})$.
We will assume that $d\equiv\ 1 (\mod\ 3)$, which ensures that
$\qq_0=(3)$ is a {\it prime} ideal in $D= \OOO_d$. Then $D/\qq_0
\cong \FF_9$ and
$$G/\pm G(\qq_0) \cong PSL_2(\FF_9)\cong A_6.$$
Hence there exists a (non-normal) subgroup $H$ of $G$ for which
$$H/\pm G(\qq_0) \cong A_5.$$
So $H$ is a congruence subgroup of level $\qq_0$ and index $6$.
\par
\noindent Let $g \in G$. Then $\qq_0\subseteq\bb(H,g)\subseteq
D$.  Suppose that $\bb(H,g)=D$. Then $T(x) \in H^g$, for all $x
\in D$, which implies that $A_5$ has a subgroup of order $9$.
Hence $\bb(H,g) \neq D$. On the other hand, by Theorem 3.1,
 $$|D:\bb(H,g)|\leq|G:H|=6.$$
It follows that $|D:\bb(H,g)|=3$.
\par
Now suppose moreover that $d\neq 1$; then $D^*=\{\pm 1\}$.
With the above notation therefore, there exists $\lambda
\in\Lambda$ such that $\bb(H,g)=\pm\bb(H,g_{\lambda})$. By Theorem
3.1 we deduce that
$$|G:H|=\sum_{\lambda\in\Lambda}|D:\bb(H,g_\lambda)|=6.$$
Hence there are precisely two $(H,G_{\infty})$-double cosets in
$G$. Then$$\bb(H,I_2)=\{0,b,-b\}+\qq_0,$$ for some $b \in D
\setminus \qq_0$. For a representative of the other double coset
 we fix $\alpha \in D$ which maps onto a primitive root of
 $\FF_9 \cong D/\qq_0$. By Lemma 3.4 there exists $g_0 \in G$ for which
 $$\bb(H,g_0)=\{0,\alpha^2b,-\alpha^2b\}+\qq_0 \neq \bb(H,I_2).$$
 Then $g_0$ can represent the other double coset. We conclude that
 $$ \{\bb(H,g): g \in G\}=\{\pm\{0,b,-b\}+\qq_0,\;\pm\{0,\alpha^2b,-\alpha^2b\}+\qq_0\}.$$
Clearly this set has neither a minimum nor maximum member under
set-theoretic containment.
\\ \\
{\bf Example 3.6.}
For our second example of this type we return to Example 2.13. Let
$\pp$ be a prime ideal in $k[t]$, where $k$ is a field. We recall that
$$ G_0(\pp)=\langle T(\alpha,\beta)\;(\alpha \in k^*,\beta\in
k),\;SL_2(k[t],\pp)\rangle.$$\\
\noindent From Example 2.13 it follows that, for all $g \in G_0$,
$$\bb(G_0(\pp),g)=a^2k+\pp=\bb_0(a),\ \mathrm{say},$$
for some $a \in k[t]$. It is clear that $\bb_0(a)=\bb_0(b)$, if $a
\equiv b\;(\mod \pp)$, and that $\bb_0(a)=\pp$ if and only if $ a
\in \pp$. We now restrict our attention to those $\bb_0(a)$ for
which $ a \notin \pp$.\\ \\
\noindent Let $K=k[t]/\pp$. For each pair $x,y \in K^*$ we write
$$x\rho y\;\mathrm{if\; and\; only\; if}\;(xy^{-1})^2 \in k^*.$$
Then $\rho$ is an equivalence relation on $K^*$. Let
$\{a_{\omega}:\omega \in \Omega\}$ be a subset of $k[t]$ which maps
bijectively onto a complete set of representatives for the
$\rho$-classes.  It follows that
$$\{\bb(G_0(\pp),g): g \in SL_2(k[t])\}=\{\bb_o(a_{\omega}): \omega
\in \Omega\}\cup\{\pp\}.$$ This set always has a minimal member,
i.e. $\pp$. However it has a maximal member only when it reduces to
$\{k+\pp,\pp\}$. This happens only when $x^2 \in k^*$, for all $x
\in K^*$. It is easy to find examples of $K$ without this property
when $k$ is not algebraically closed.

\subsection*{4. Quasi-level}

\noindent{\bf Definition.} Let $H$ be a subgroup of $G$. We define
the {\bf quasi-level} of $H$ as
$$ql(H)=\displaystyle{\bigcap_{g \in G}}\bb(H,g).$$
\\
Since $\bb(H,g)=\bb(H^g,I_2)$, we see that
$$ql(H)=\displaystyle{\bigcap_{g \in G}}\bb(H^g,I_2)=\bb(N_H,I_2)$$
where $N_H$ is the core of $H$ in $G$.
\par
However, Example 3.2 shows that in contrast to Lemma 2.3 the 
intersection over $g\in G$ can in general not be replaced by 
the intersection over a system of representatives of the cusps 
(even when $H$ is a congruence subgroup).
\par
The second equality shows that our definition of quasi-level coincides
with the one we gave in [MSch2]. Actually $ql(H)=\bb(N_H,g)$ for any
$g\in G$. But Example 3.5 shows that in contrast to Corollary 2.11 there
is in general no $g\in G$ for which $ql(H)=\bb(H,g)$ (even when $H$ is
a congruence subgroup).
\par
\noindent We summarize the basic properties:
\\ \\
{\bf Lemma 4.1.} \it \begin{itemize}
\item[(i)] $ql(H)=\{d \in D:T(d) \in N_H\}.$
\item[(ii)] $ql(H)$ is an additive subgroup of $D$ with the
property that
$$ d \in ql(H),\;u \in D^*\Longrightarrow u^2d \in ql(H).$$
\item[(iii)] $ql(H)\supseteq l(H).$
\item[(iv)] $ql(H)=ql(H^g)=ql(N_H)=\bb(N_H,g)$ for all $g \in G$.
\item[(v)] $l(H)=l(H^g)=l(N_H)=\cc(N_H,g)$ for all $g \in G$.
\item[(vi)] $l(H)$ is the largest $D$-ideal contained in $ql(H)$.
\item[(vii)] $|D:ql(H)|\le |SL_2(D):N_H|$.
\end{itemize}

\noindent \rm For congruence subgroups we can combine Lemmas 3.4 and
4.1 to obtain the following extension of Lemma 4.1(ii).
\\ \\
{\bf Lemma 4.2.} \it Let $H$ be a congruence subgroup and
let $\alpha \in D$ be invertible modulo $l(H)$. Then $$\alpha^2ql(H)
\subseteq ql(H).$$

\noindent \rm We will also show that when $H$ is a congruence
subgroup the inequality in Lemma 4.1(iii) becomes an equality (in
``most" cases). For this purpose we require a number of
preliminaries.
\\ \\
{\bf Lemma 4.3.} \it Let $H$ be a congruence subgroup of $G$ and let
$l(H)=\qq=\qq_1\qq_2$, where $\qq_1+\qq_2=D$. Then
$$ l((H\cap G(\qq_1))\cdot G(\qq_2))=\qq_2.$$
\\
{\bf Proof.} \rm Let the required level be $\qq_2'$. Then $\qq_2'
\supseteq \qq_2$. \noindent Now
$$ G(\qq_2') \leq (H\cap G(\qq_1))\cdot G(\qq_2)$$
and so
$$G(\qq_1\qq_2')=G(\qq_1)\cap G(\qq_2')\leq (H\cap G(\qq_1))\cdot
G(\qq)\subseteq H.$$
\noindent Hence $\qq\supseteq\qq_1\qq_2'$ by Corollary
1.4. The result follows.
\hfill $\Box$
\\ \\
{\bf Lemma 4.4.} \it Let $N$ be a normal congruence
subgroup of level $\qq=\qq_1\qq_2$, where $\qq_1+\qq_2=D$. Let
$N_0=(N\cap G(\qq_1))\cdot G(\qq_2)$ and $\overline{N}=N \cdot
G(\qq_2)$. Then

$$\overline{N}/N_0 \;is\; a\; central\; subgroup\; of\; G/N_0.$$

\noindent {\bf Proof.} \rm Now $G=G(\qq_1)\cdot G(\qq_2)$, by
Corollary 1.3. It follows that
$$[G,\overline{N}]=[G(\qq_1)\cdot G(\qq_2),N\cdot G(\qq_2)]\leq [G(\qq_1),N]\cdot G(\qq_2) \leq N_0.$$
The result follows. \hfill $\Box$
\\ \\
\noindent Our next lemma is almost certainly well-known. In the
absence of a reference we provide a proof.
\\ \\
{\bf Lemma 4.5.} \it Let $L$ be a local ring for which $2
\in L^*$. Then $PSL_2(L)$ has trivial centre.
\\ \\
\noindent {\bf Proof.} \rm  Note that, since $2 \in L^*$, the only
involutions in $L^*$ are $\pm 1$. Let
$$g=\left[\begin{array}{cc}\alpha&\beta\\
\gamma&\delta\end{array}\right] \in SL_2(L).$$\\
\noindent Then, if $g$ maps into the centre of $PSL_2(L)$, it
follows that, for all $x \in L$,
$$ gT(x)=\lambda T(x)g \ \ \mathrm{and}\ \ gS(x)=\mu S(x)g,$$
where $\lambda^2=\mu^2=1$. \noindent If $ \gamma \in L^*$, then
$\lambda =1$ and so $\gamma x=0$. Thus $\gamma \notin L^*$ and
similarly $\beta \notin L^*$. We deduce that $\alpha,\delta \in
L^*$. \noindent From the first of the above equations it follows
that $\alpha=\lambda(\alpha+x\gamma)$ and hence that
$2x\alpha\gamma+x^2\gamma^2=0$. The latter equation holds for $x=\pm
1$  and so $4\alpha\gamma=0$. From the above $\gamma=0$ and
similarly $\beta=0$. It then follows from the above that
$\alpha=\delta$, i.e. $g=\pm I_2$.
\hfill $\Box$
\\ \\
\noindent Before coming to our next principal result we make
another definition.
\\ \\
{\bf Definition.}
For a subgroup $H$ of $G=SL_2(D)$ we define $o(H)$ as the ideal of $D$
generated by all elements $a-d,\ b,\ c$ with ${a\ b\choose c\ d}\in H$.
Somewhat unfortunately, $o(H)$ is sometimes called the {\bf order} of $H$.
Every matrix in $H$ is congruent modulo $o(H)$ to a scalar matrix $xI_2$
for some $x\in D$, and $o(H)$ is the smallest ideal of $D$ with this
property. Obviously
$$ql(H)\subseteq o(H).$$
Conversely, for each $D$-ideal $\qq$ we define
$$ Z(\qq)=\{ X \in G: X \equiv xI_2\;(\mod \qq)\;
\mathrm{for\; some}\; x \in D\}.$$
Then
$$H\leq Z(\qq)\Leftrightarrow o(H)\leq\qq.$$
So it is clear that if
$$ G(\qq) \leq H \leq Z(\qq),$$
then $$ ql(H)=l(H)=\qq.$$
\\
{\bf Theorem 4.6.} \it Let $H$ be a congruence subgroup of $SL_2(D)$
such that $l(H)=\qq$ satisfies Condition L (from the Introduction).
Then
$$ql(H)=l(H),$$
equivalently, the quasi-level is actually an ideal.
\\ \\
\noindent {\bf Proof.} \rm By Lemma 4.1(iv), (v) we may assume that
$H=N_H$ (i.e. $H\unlhd G$). From the above it is sufficient to prove
that $ N_H \leq Z(\qq)$, and hence that $N_H \leq Z(\pp^{\alpha})$,
where $\pp$ is any prime ideal
for which $\alpha=\mathrm{ord}_{\pp}(\qq) >0$.\\
\noindent Let $L$ denote the local ring $D/\pp^{\alpha}$ and let
$$\pi:G \longrightarrow SL_2(L),$$
denote the natural map. Now $SL_2(L)$ is generated by elementary
matrices by [K, Theorem 1] and so (again by [K, Theorem 1])
$$\pi(G(\rr))=E_2(L,\overline{\rr})=SL_2(L,\overline{\rr}),$$
for all $D$-ideals $\rr$, where $\overline{\rr}$ is the image of
$\rr$ in $L$. (In particular $\pi$ is an {\it
epimorphism}.)\\
\noindent Let $\qq=\qq'\pp^{\alpha}$. Now suppose that $\pi(N_H\cap
G(\qq'))$ is not central in $SL_2(L)$. Then by [K, Theorem 3] the
hypotheses on $\qq$ and the above ensure that
$$ G(\pp^{\beta})\leq (N_H\cap G(\qq'))\cdot G(\pp^{\alpha}).$$
for some $\beta < \alpha$ which contradicts Lemma 4.3 (with
$\qq_2=\pp^{\alpha}$ and $\qq_1=\qq'$). It follows that $$ \pi(N_H
\cap G(\qq')) \leq \{\pm I_2\}.$$ The map $\pi$ extends to an
epimorphism $$\overline{\pi}: G \longrightarrow PSL_2(L).$$ By Lemma
4.4 $\;\overline{\pi}(N_H)$ is central in $PSL_2(L)$. We now apply
Lemma 4.5 to conclude that
$$\pi(N_H) \leq \{\pm I_2\}.$$ \hfill $\Box$
\\ \\
For normal subgroups we can reformulate Theorem 4.6 as follows.
\\ \\
{\bf Corollary 4.7.} \it
Let $N$ be a normal congruence subgroup of $G$. If the level of $N$
satisfies Condition L, then
$$\bb(N,g)=\cc(N,g),$$
that is, the quasi-amplitudes of $N$ are actually the cusp amplitudes.
\rm
\\ \\
{\bf Proof.} \rm
Since $N$ is normal any quasi-amplitude is equal to the quasi-level
and any cusp amplitude is equal to the level.
\hfill $\Box$
\\ \\
{\bf Remark 4.8.}
A word of warning is in order here. If $ql(H)=l(H)$ for a (non-normal)
congruence subgroup $H$ and $H$ has only one cusp, this does {\it not}
imply $\bb(H,g)=\cc(H,g)$, not even if $l(H)$ satisfies Condition L.
See Example 3.2. Ultimately the problem is caused by diagonal matrices
that are not central.
\\ \\
\noindent McQuillan [Mc, Theorem 1] has proved, for the special case
$D=\ZZ$, that, if $N$ is a normal congruence subgroup of $G$ of
level $\qq$, then $ N \leq Z(\qq)$, using a similar approach. We now
provide a pair of examples to show that both
restrictions in Theorem 4.6 are necessary.
\\ \\
{\bf Example 4.9.}
Our first example [M1, Example 2.3] shows that Theorem 4.6 can fail
when $\qq$ is not prime to $2$. Let $\pp$ be a prime $D$-ideal for
which $2 \in \pp^2$. We recall from Lemma 1.1 that
$$ G(\pp^2)/G(\pp^4) \cong \{(a,b,c): a,b,c \in \pp^2/\pp^4\},$$
\noindent where the latter group is additive. Let $\Lambda = \{
t^2+\pp^4: t \in \pp\}$. We define a subgroup $K$, where $G(\pp^4)
\leq K \leq G(\pp^2)$, by
$$ K/G(\pp^4)= \{(a,b,c): b,c \in \Lambda\}.$$
\noindent Since $2\pp^2 \subseteq \pp^4$, it is easily verified that $K$
is a subgroup of $G$, normalized by $S(x),T(x)$, for all $x \in D$.
Since $K$ is also normalized by $G(\pp^4)$, it follows from Theorem
1.2 that $K$ is {\it normal} in $G$.
\par
Clearly $ql(K)=\{t^2+q:t \in \pp,q \in \pp^4\}$. Now
suppose that $l(K) \neq \pp^4$. Then $G(\pp^3) \leq K $. Let $h$ be
a generator of $\pp^3\;(\mod \pp^4)$. It follows that there exists
$k \in D$ such that $h \equiv k^2\; (\mod \pp^4)$. We conclude that
$l(K)=\pp^4$.
\\ \\
Explicit examples that satisfy the requirement $2\in\pp^2$ are,
among others, $D=\ZZ[\sqrt{-2}]$ with $\pp=(\sqrt{-2})$, or to take
a local example, $D=\ZZ_2[\sqrt{2}]$ with $\pp=(\sqrt{2})$. More
generally, this example actually works for any Dedekind domain of
characteristic $2$ with any nonzero prime ideal $\pp$ because then
trivially $2=0\in\pp$.
\\ \\
{\bf Example 4.10.}
Suppose that $\qq=\mm_1\mm_2$, where $\mm_1+\mm_2=D$ and
$|D/\mm_i|=3\;(i=1,2)$. (Consider, for example, $D=\ZZ[\sqrt{-2}]$
with $\mm_1=(1+\sqrt{-2})$ and $\mm_2=(1-\sqrt{-2})$ or $D=\FF_3[t]$
with $\mm_1=(t)$ and $\mm_2=(t+1)$.) Then, by Corollary 1.3(ii),
$$G/G(\qq) \cong SL_2(\FF_3) \times SL_2(\FF_3).$$
\noindent From the well-known structure of $SL_2(\FF_3)$ it follows
that there exists a {\it normal} subgroup, $N$, of $G$, containing
$G(\qq)$, such that
$$ |G:N|=9\;\;\mathrm{and}\;\;|N:G(\qq)|=64.$$
\noindent Now let
$$ M=\langle T(1),N \rangle.$$
\noindent Since $9=3^2$, $M \unlhd G$ and $|G:M|=3$. Obviously,
$ql(M)$ contains $1$. If $ 1 \in l(M)$, then $l(M)=D$, in which case
$M=G(D)=G$. Thus $ql(M)\neq l(M)$.
\par
In particular, there exists a normal congruence subgroup of index
$3$ in $SL_2(\FF_3[t])$ that has level $t(t+1)$.
\\ \\
{\bf Remarks 4.11.}
\begin{itemize}
\item[a)] If the level of a congruence subgroup $H$ is a prime ideal
$\pp$, then $ql(H)=l(H)$. This follows immediately from the
simplicity of the group $PSL_2(D/\pp)$, when $|D/\pp|> 3$. The cases
for which $|D/\pp|\le 3$ can be checked directly.
\item[b)] If $D$ is any arithmetic Dedekind domain, the quasi-level
of a congruence subgroup is not ``too far from" its level.
For a normal congruence subgroup $N$ the relation between $l(N)$ and
$o(N)$ is described in [M2, Theorems 3.6, 3.10 and 3.14]. See also
the end of Section 3 of [M2].
\par
For example, for a finite index subgroup $H$ of $SL_2(\ZZ[\sqrt{11}])$
we obtain
$$4ql(N_H)\subseteq 4o(N_H)\subseteq l(N_H)\subseteq ql(N_H)\subseteq o(N_H)$$
from [M2, Theorem 3.6] since $2$ is ramified and $3$ is inert in
$\ZZ[\sqrt{11}]$. Actually even $4ql(H)\subsetneq l(H)$
since $ql(N_H)=o(N_H)$ would mean that $ql(H)$ is an ideal and hence
equal to $l(H)$.
\end{itemize}

\noindent Our final result demonstrates that for a non-congruence
subgroup there is in general almost no connection between its
quasi-level and level (in contrast with Theorem 4.6). We note that
there is {\it no} proper normal subgroup of $SL_2(k[t])$ whose
quasi-level is $k[t]$. (Since $k[t]$ is a Euclidean ring,
$SL_2(k[t])$ is generated by $T(r),S(r)$, where $r \in k[t]$.)
\\ \\
{\bf Theorem 4.12.} \it
Let $k$ be any field and let $f\in k[t]$ with $\deg(f)\ge 2$.
Suppose that $f(0) \neq 0$ and, further, that $f'(0) \neq 0$, when
$\deg(f)=2$. Then there exists a normal non-congruence subgroup $N$
of $SL_2(k[t])$ of level $(f)$ with the following properties.
\begin{itemize}
\item[{(i)}] $N\cdot SL_2(k)=SL_2(k[t])$.
\item[{(ii)}] $l(N)=(f)$.
\item[{(iii)}] $ql(N)$ has $k$-codimension $1$ in $k[t]$.
\end{itemize}
\noindent {\bf Proof.} \rm We note that, by hypothesis, $t \nmid f$.
We define the $k$-subspace
$$Q=(f)\oplus kt\oplus kt^2\oplus\cdots\oplus kt^{d-1}$$
where $d=\deg(f)$. Let $N=\Delta(Q)$ be the {\it normal} subgroup of
$SL_2(k[t])$ generated by all $T(q)$, where $q \in Q$. Since
$SL_2(k[t])$ is generated by all $T(r),S(r)$, where $r \in k[t]$,
part (i) follows. In addition
$$ql(\Delta(Q))=Q,$$
by [M3, Theorem 3.8]. Part (iii) follows. \noindent Suppose that
$l(\Delta(Q)) \neq (f)$. Then
$$ (f) \subseteq (h) \subseteq Q,$$ for some polynomial divisor $h$
of $f$, with $\deg(h) < \deg(f)$. Then by the definition of $Q$,
$h$, and hence $f$, must be divisible
by $t$. Part (ii) follows.
\par
Finally, suppose that $N$ is a congruence subgroup.
Since, by hypothesis $t$ is prime to $l(N)=(f)$, $$ t^2Q \subseteq
Q,$$ by Lemma 4.2. If $\deg(f)>2$, then $t^d \in Q$. If $\deg(f)=2$,
then $t^3 \in Q$. Now for this case $tf \in Q$ and so $t^2 \in Q$,
by the extra hypothesis. In either case $(t) \subseteq Q$, which
implies that $Q=k[t]$, a contradiction. The proof is complete.
\hfill $\Box$
\\ \\
{\bf Remarks 4.13.}
\begin{itemize}
\item[a)] Obviously the group $N$ in Theorem 4.12 shows that Lemmas
3.4 and 4.2 do not hold in general for non-congruence subgroups,
even if they are normal.
\item[b)] The restriction on the degree of $f$ in Theorem 4.12 is
necessary. It is well-known that, if $\deg(f) \leq 1$, then every
subgroup of $SL_2(k[t])$ of level $(f)$ is a congruence subgroup.
\item[c)] Several versions of Theorem 4.12 are already known for
subgroups of level zero. (See Section 4 of [MSch2].)
\end{itemize}

\subsection*{5. Level and index}

\noindent We now make use of quasi-amplitudes to extend the {\it
index/level} inequality $(*)$ for finite index subgroups of
$SL_2(\ZZ)$ (in Example 2.12) to other arithmetic Dedekind domains.
\\ \\
First, let $K$ be an algebraic number field of degree $d$ over $\QQ$.
We recall the definition of $\OOO_S$ from Section 1, where $S$ is a
suitable, not necessarily finite, set of places of $K$.
\\ \\
{\bf Lemma 5.1.} \it
Let $L$ be a subgroup of $\OOO_S$ such that $\OOO_S/L$ has exponent
$e$. Then $L$ contains the $\OOO_S$-ideal generated by $e$.
\par
In particular, every finite index subgroup $L$ of $\OOO_S$ contains
an $\OOO_S$-ideal $\aa$ such that
$$|\OOO_S/\aa|\leq |\OOO_S:L|^d.$$
\rm\\
{\bf Proof.}
The first claim is obvious since $ea\in L$ for every $a\in\OOO_S$.
\par
\noindent For the second claim let $|\OOO_S:L|=n$. Then $L$ contains
$n\OOO_S$. Denote the ring of integers of $K$ by $\OOO_K$. If
$n\OOO_K=\prod\pp_i^{e_i}$ is the decomposition of $n\OOO_K$ into a
product of prime ideals of $\OOO_K$, then
$$n\OOO_S =\prod_{\pp_i\not\in S}\pp_i^{e_i}\OOO_S.$$
So
$$|\OOO_S/n\OOO_S|=\prod_{\pp_i\not\in S}|\OOO_S/\pp_i\OOO_S|^{e_i}=
\prod_{\pp_i\not\in S}|\OOO_K/\pp_i\OOO_K|^{e_i}$$
$$\le\prod|\OOO_K/\pp_i\OOO_K|^{e_i}=|\OOO_K/n\OOO_K|=n^d.$$
\hfill $\Box$
\\ \\
{\bf Theorem 5.2.} \it
Let $K$ be an algebraic number field with $[K:\QQ]=d$ and let
$\OOO_S$ be the ring of $S$-integers of $K$ (with $|S|$ not necessarily
finite). If $H$ is a congruence subgroup of $SL_2(\OOO_S)$, then
$$|\OOO_S/l(H)|\leq |SL_2(\OOO_S):H|^d.$$
\rm
\\
\noindent {\bf Proof.} \rm
By Corollary 2.11 there exists $g\in SL_2(\OOO_S)$ with
$l(H)=\cc(H,g)$. Let $L$ be the quasi-amplitude $\bb(H,g)$. Then
$|\OOO_S:L|\leq |SL_2(\OOO_S):H|$, and the Theorem follows
from Lemma 5.1.
\hfill $\Box$
\\ \\
{\bf Remarks 5.3.}
\begin{itemize}
\item[a)] Theorem 5.2 also holds for normal non-congruence subgroups
of $SL_2(\OOO_S)$ since then every cusp amplitude is equal to the level.
\item[b)] If $|S|>1$ and $S$ contains at least one real or non-archimedean
place, or if $|S|=\infty$, then by [Se1] resp. Theorem 1.9 the inequality
in Theorem 5.2 holds for {\it all} subgroups of finite index in
$SL_2(\OOO_S)$.
\item[c)] Lubotzky [Lu, (1.6) Lemma] has given a version of Theorem 5.2
for more general algebraic groups using different methods.
\end{itemize}
\noindent Next we exhibit some examples for which the inequality in
Theorem 5.2 is sharp.
\\ \\
{\bf Example 5.4.}
Let $K$ be a number field of degree
$d$ over $\QQ$. Let $D=\OOO_K$ denote the ring of integers of $K$
(i.e. the ring of $S$-integers of $K$, where $S$ consists precisely
of the archimedean places of $K$). We suppose that the ideal
$\qq_0=(2)=2D$ splits into the product of $d$ distinct prime
ideals in $D$. Let $G=SL_2(D)$. Then, using Corollary 1.3,
$$G/G(\qq_0)\cong P_1 \times\cdots\times P_d,$$
where for $1 \leq i \leq d$,$$ P_i \cong SL_2(\FF_2) \cong S_3,$$
Let $\pp$ be any prime ideal dividing $\qq_0$ and let
$\qq_0=\pp\pp'$. Then, again by Corollary 1.3, under the first
isomorphism above $$ G(\pp')/G(\qq_0) \cong P_j,$$ for some $j$.
\noindent For each $i$, let $N_i$ be the normal subgroup of $P_i$ of
order $3$. There exists an epimorphism $$\theta:
P_1\times\cdots\times P_d \longrightarrow S_2,$$ such that $N=
\mathrm{ker}\;\theta$ contains $N_1 \times \cdots \times N_d$ but
{\it not} any $P_i$. Let $M/G(\qq_0)$ be the inverse image of
$N$ in $G/G(\qq_0)$. If $l(M)\neq \qq_0$, then, by Corollary
1.4, $M$ contains $G(\pp')$, for some prime divisor $\pp$ of
$\qq_0$, which contradicts the above. Hence $l(M)=\qq_0$. We
conclude that $M$ is a normal congruence subgroup with
$$ |D/l(M)|=|G:M|^d=2^d.$$
By the way, $|G:M|=2$ also shows
$$[D:ql(M)|=|D:\bb(M,g)|=2.$$
\\
{\bf Remark 5.5.}
Actually, for each $d$ there are infinitely many number fields
$K$ of degree $d$ that satisfy the condition in Example 5.4.
This can be seen as follows:
\par
\noindent By Dirichlet's Theorem on primes in arithmetic
progressions there are infinitely many odd primes that are congruent
to $1$ modulo $d$. Pick two of these, say $p_1$ and $p_2$ with
$p_i=m_i d+1$. Then $2^{m_id}\equiv 1\ \mod\ p_i$ and hence
$2^{m_1m_2d}\equiv 1\ \mod\ p_1p_2$. So by the decomposition law in
cyclotomic fields the inertia degree of $(2)$ in the $p_1p_2$-th
cyclotomic field $\QQ(\zeta_{p_1p_2})$ divides $m_1m_2d$. Since
$\QQ(\zeta_{p_1p_2})$ has degree $m_1m_2d^2$ over $\QQ$, the
decomposition field of $(2)$ in $\QQ(\zeta_{p_1p_2})$ is an abelian
extension of $\QQ$ whose degree is divisible by $d$. Thus it
contains a subfield $K$ of degree $d$ over $\QQ$ with the desired
property.
\par
Alternatively we could argue as follows:
$Gal(\QQ(\zeta_{p_1p_2})/\QQ)\cong\ZZ/m_1 d\ZZ\oplus\ZZ/m_2 d\ZZ$ and
the Frobenius at $2$ generates a cyclic subgroup; so its fixed field
is an abelian extension of $\QQ$ whose degree is divisible by $d$.
\par
Taking two other primes $p_1$, $p_2$ one obtains a different field
$K$ since $\QQ(\zeta_{p_1p_2})$ and hence $K$ is unramified outside
$p_1p_2$.
\\ \\
By an analogous proof for each $d$ there are infinitely many number
fields $K$ of degree $d$ such that $(3)$ splits completely in $\OOO_K$.
For these one can similarly construct a normal congruence subgroup of
index $3$ in $SL_2(\OOO_K)$ with level $3\OOO_K$.
(Compare Example 4.10 for a special case.)
\\ \\
Thus it looks like Theorem 5.2 is optimal. Also, it is not possible
to prove a function field analogue of Theorem 5.2 along the same lines,
as there cannot be a function field analogue of the inequality in
Lemma 5.1. Even in $\FF_q[t]$ one can construct additive subgroups
$L$ of index $q$ such that the biggest ideal $\aa$ contained in $L$
is the zero ideal or any prescribed nontrivial ideal.
(This is actually the key point for many constructions in [MSch2].
See also Theorem 4.12.)
\par
However, somewhat surprisingly, under the condition that the level is
prime to certain ideals one can give a relation between the level and
the index (if finite) of a congruence subgroup of $SL_2(D)$ that is
valid for {\it any} Dedekind domain $D$.
\\ \\
{\bf Theorem 5.6.} \it
Let $D$ be any Dedekind domain and $N$ a normal congruence subgroup
of $SL_2(D)$. If the level of $N$ satisfies Condition L, then
$$|D/l(N)|\ \hbox{\it divides}\ |SL_2(D):N|.$$
\rm
\\
\noindent {\bf Proof.} \rm
If $N$ is normal, the level is equal to any cusp amplitude. Under
Condition L this cusp amplitude is equal to the quasi-amplitude
by Corollary 4.7. So the result follows from Theorem 3.1.
\hfill $\Box$
\\ \\
Examples 3.2 and 3.5 show that even for ``relatively simple'' rings
like $\FF_9[t]$ or $\ZZ[\sqrt{-13}]$ we cannot expect
$$|D/l(H)|\le |SL_2(D):H|$$
for non-normal congruence subgroups $H$ that satisfy Condition L.
But of course we have
\\ \\
{\bf Corollary 5.7.} \it
Let $D$ be any Dedekind domain and $H$
a congruence subgroup of $SL_2(D)$ of index $n$. If the level of $H$
satisfies Condition L, then
$$|D/l(H)|\ \hbox{\it divides}\ n!.$$
\rm
\\
\noindent {\bf Proof.} \rm
$H$ has the same level as its core in $SL_2(D)$, and the index of this
core divides $n!$.
\hfill $\Box$
\\ \\
Corollary 5.7 is stronger than Theorem 5.2 if $n$ is sufficiently small
compared to $d$.
\par
Examples 4.10 and 5.4 show that Condition L cannot simply be dropped in
Theorem 5.6 and Corollary 5.7.
But for arithmetic Dedekind domains $D$ one can prove somewhat weaker results
than Theorem 5.6 and Corollary 5.7 for congruence subgroups that do not
satisfy Condition L. The key is that by results of [M2] one can still
control the relation between $l(N)$ and $ql(N)$. See our Remark 4.11 b).
We content ourselves with the function field case. Since the case of
characteristic $p>3$ is already fully covered by Theorem 5.6, we only
have to deal with characteristic $2$ and $3$.
\\ \\
{\bf Theorem 5.8.} \it
Let $K$ be a global function field with constant field $\FF_q$ and let
$D=\OOO_S$ where $|S|$ is finite. Let $N$ be a normal congruence subgroup
of $SL_2(D)$.
\begin{itemize}
\item[a)] If $char(\FF_q)=3$, then
$$|D/l(N)|\ \hbox{\it divides}\ 3^r |SL_2(D):N|$$
where $r$ is the number of prime ideals $\pp$ of $D$ with $|D/\pp|=3$.
\item[b)] If $char(\FF_q)=2$, then
$$|D/l(N)|\ \hbox{\it divides}\ 4^s |SL_2(D):N|^2$$
where $s$ is the number of prime ideals $\pp$ of $D$ with $|D/\pp|=2$.
\end{itemize}
\rm

\noindent {\bf Proof.} \rm
We only give the proof for characteristic $2$ as they are almost the same.
\par
By [M2, Theorem 3.14] (compare also the end of Section 3 in [M2]) we have
$\ss^2(o(N))^2\leq l(N)$ where $\ss$ is the product of all prime ideals
$\pp$ in $D$ with $|D/\pp|=2$. Hence $|D/l(N)|$ divides $4^s |D/o(N)|^2$,
which divides $4^s |D/ql(N)|^2$ since $ql(N)$ is an $\FF_q$-subspace of
$o(N)$. As $N$ is normal, we have $ql(N)=\bb(N,g)$, and the claim follows
from Theorem 3.1.
\hfill $\Box$
\\ \\
{\bf Remarks 5.9.}
\begin{itemize}
\item[a)]
If $K$ is a number field and $|S|$ is finite, then using [M2, Theorem 3.6]
and the end of Section 3 of [M2], for a normal congruence subgroup $N$ of
$SL_2(\OOO_S)$ we obtain in the worst case
$$|\OOO_S /l(N)|\ \hbox{\rm divides}\ 12^d |SL_2(\OOO_S) :N|$$
where $d=|K:\QQ|$.
\item[b)]
Theorems 5.6 and 5.8 can be used to improve [MSch2, Proposition 4.5 b)] in
the sense that the condition that $\mm$ is a maximal ideal is not really
needed there. This condition was used to show that the constructed groups
of level $\mm$ are non-congruence subgroups. But if $\deg(\mm)$ is big enough,
then Theorem 5.6 (resp. Theorem 5.8) guarantees the non-congruence property.
\end{itemize}
\bigskip
\noindent {\bf Acknowledgements.} 
Parts of this paper were developed during research visits of the second
author at Glasgow University. The hospitality of their Mathematics Department
is gratefully acknowledged. The final version was written while the second 
author was holding a visiting position at the National Center for Theoretical
Sciences (NCTS) in Hsinchu, Taiwan. We are also grateful to Peter V\'amos for 
pointing out the reference [Go] to us.
\\

\subsection*{\hspace*{10.5em} References}

\begin{itemize}
\item[{[AS]}] A. O. L. Atkin and H. P. F. Swinnerton-Dyer:
Modular forms on noncongruence subgroups, \it
Combinatorics (Proc. Sympos. Pure Math., Vol. XIX, Univ. California,
Los Angeles, Calif., 1968), \rm pp.1-25.
\it Amer. Math. Soc., Providence, R.I., $1971$.\rm
\item[{[B]}] H. Bass: \it Algebraic $K$-theory,
\rm W. A. Benjamin Inc., New York-Amsterdam, 1968.\rm
\item[{[Ge]}] E.-U. Gekeler: \it Drinfeld Modular Curves, \rm
Lecture Notes in Mathematics, vol. 1231, Springer, Berlin Heidelberg
New York, 1986.
\item[{[Go]}] O. Goldman: On a special class of Dedekind domains,
\it Topology \bf 3 suppl. 1 \rm (1964), 113-118.
\item[{[K]}] W. Klingenberg: Lineare Gruppen $\ddot{\mathrm{u}}$ber lokalen
Ringen, \it Amer. J. Math. \bf 83 \rm (1961), 137-153.
\item[{[La]}] H. Larcher: The cusp amplitudes of the congruence subgroups
of the classical modular group,
\it Illinois J. Math. \bf 26 \rm (1982), 164-172.
\item[{[Lu]}] A. Lubotzky: Subgroup growth and congruence subgroups,
\it Invent. Math. \bf 119 \rm(1995), 267-295.
\item[{[Mc]}] D. L. McQuillan: Classification of normal congruence
subgroups of the modular group, \it Amer. J. Math. \bf 87 \rm
(1965), 285-296.
\item[{[M1]}] A. W. Mason: Standard subgroups of $GL_2(A)$, \it
Proc. Edinburgh Math. Soc. \bf 30 \rm (1987), 341-349
\item[{[M2]}] A. W. Mason: The order and the level of a subgroup of
$GL_2$ over a Dedekind ring of arithmetic type, \it Proc. Roy. Soc.
Edinburgh Sect. A \bf 119 \rm (1991), 191-212.
\item[{[M3]}] A. W. Mason: Normal subgroups of $SL_2(k[t])$ with or
without free quotients, \it J. Algebra \bf 150 \rm (1992), 281-295.
\item[{[MSch1]}] A. W. Mason and A. Schweizer: The minimum index
of a non-congruence subgroup of $SL_2$ over an arithmetic domain.
II: The rank zero cases, \it J. London Math. Soc. \bf 71 \rm
(2005), 53-68
\item[{[MSch2]}] A. W. Mason and A. Schweizer: Non-standard automorphisms
and non-cong\-ruence subgroups of $SL_2$ over Dedekind domains contained
in function fields, \it J. Pure Appl. Algebra \bf 205 \rm (2006), 189-209
\item[{[MSt]}] A. W. Mason and W. W. Stothers:
On subgroups of $GL(n,A)$ which are generated by commutators,
\it Invent. Math. \bf 23 \rm (1974), 327-346.
\item[{[Se1]}] J.-P. Serre: Le probl\`eme des groupes de congruence
pour $SL_2$, \it Ann. of Math. \bf 92 \rm (1970), 489-527.
\item[{[Se2]}] J.-P. Serre: \it Trees, \rm Springer, Berlin
Heidelberg New York, 1980
\item[{[St]}] W. W. Stothers: Level and index in the modular group,
\it Proc. Roy. Soc. Edinburgh Sect. A \bf 99 \rm(1984), 115-126.
\item[{[Ta]}] S. Takahashi: The fundamental domain of the
tree of $GL(2)$ over the function field of an elliptic curve,
\it Duke Math. J. \bf 72 \rm (1993), 85-97
\item[{[W]}] K. Wohlfahrt: An extension of F. Klein's level concept,
\it Illinois J. Math. \bf 8 \rm(1964), 529-535.
\end{itemize}

\end{document}